\documentclass[a4paper,12pt]{article}

\usepackage[top=1in,left=1in,right=1in,bottom=1.5in]{geometry}

\usepackage[T1]{fontenc}
\usepackage{amsmath,amssymb,amsfonts,amsthm,mathrsfs,eucal,dsfont,verbatim}

\evensidemargin
-1cm

\usepackage{hyperref}

\usepackage{color}

\newcommand{\normal}{\color{black}}

\theoremstyle{plain}
\newtheorem{theorem}{Theorem}[section]
\newtheorem{lemma}[theorem]{Lemma}
\newtheorem{proposition}[theorem]{Proposition}
\newtheorem{corollary}[theorem]{Corollary}
\theoremstyle{definition}
\newtheorem{remark}{Remark}[section]

\newcommand{\prt}{\partial}
\renewcommand{\star}{\circledast}

\newcommand{\real}{\mathds{R}}
\newcommand{\rn}{{\mathds{R}^n}}
\newcommand{\rd}{{\mathds{R}^d}}

\newcommand{\Ee}{\mathds{E}}
\newcommand{\I}{\mathds{1}}

\newcommand{\Pp}{\mathds{P}}

\newcommand{\Sss}{\mathbb{S}}

\newcommand{\RR}{\mathrm{Re}\,}

\newcommand{\eps}{\varepsilon}
\numberwithin{equation}{section}

\begin{document}
\title{Intrinsic compound kernel estimates for the  transition probability density of
L\'evy type processes and their applications}

\author{%
    \textsc{Victoria Knopova}%
    \thanks{  V.M.\ Glushkov Institute of Cybernetics,
            NAS of Ukraine,
            40, Acad.\ Glushkov Ave.,
            03187, Kiev, Ukraine,
            \texttt{vic\underline{ }knopova@gmx.de}}
    \textrm{\ \ and\ \ }
    \stepcounter{footnote}\stepcounter{footnote}\stepcounter{footnote}
    \stepcounter{footnote}\stepcounter{footnote}%
    \textsc{Alexei Kulik}%
    \thanks{Institute of Mathematics, NAS of Ukraine, 3, Tereshchenkivska str., 01601  Kiev, Ukraine,
    \texttt{kulik@imath.kiev.ua}}
    }

\date{}

\maketitle

\begin{abstract}
    \noindent

    Starting with an integro-differential operator $(L, C_\infty^2(\rn))$, we prove that its $C_\infty(\rn)$-closure  is the generator of a Feller process $X$, which admits a transition probability density. To construct this transition probability density, we develop a version of the parametrix method and a  verification procedure, which proves that the constructed object is the claimed one.  As a part of the construction, we prove  the intrinsic upper and lower estimates on the density.   As an application  of the constructed estimates  we state the necessary and (separately) sufficient conditions  under which a given Borel measure belongs to the  Kato   and Dynkin classes with respect to the constructed  transition probability density.


\end{abstract}

\numberwithin{equation}{section}

\section{Introduction}

 A Markov process $X=(X_t)_{t\geq 0}$ with values in  $\rn$ is called a \emph{L\'evy type process}, if its  generator $A$ is well defined on the space $C_\infty^2(\rn)$ of twice continuously differentiable functions, vanishing at $\infty$ together with their derivatives, and on this space $A$ coincides with a \emph{L\'evy type operator}
\begin{equation}\label{gener}
\begin{split}
Lf(x)&=  a(x)\cdot  \nabla f(x) + \sum_{j,k=1}^n  Q_{jk}(x)\frac{\partial^2 f(x)}{\partial x_j\partial x_k} \\
& \quad + \int_{\rn \backslash \{0\}} (f(x+u)-f(x)-\nabla f(x) \cdot u \I_{\{\|u\|\leq 1\}}\big)\mu(x,du),
\end{split}
\end{equation}
where  $a(x)\in \rn$, $Q(x)\equiv (Q_{jk}(x))_{j,k=1}^n$  is a symmetric positive semi-definite matrix, and $\mu(x,\cdot)$ is a positive Borel measure, such that $\int_{\rn\backslash \{0\}} (1\wedge \|u\|^2)\mu(x,du)<\infty $ for any $x\in \rn$.   
For an extensive survey on L\'evy type processes and L\'evy type operators we refer  to  \cite{Ja01}--\cite{Ja05}, and  \cite{BSW}; here we briefly outline the items important for the exposition below.

In the ``constant coefficients case'',  where $a,Q,\mu$ do not depend on $x$, (\ref{gener})  is just an expression for the generator of the semigroup of probability measures, which corresponds to a L\'evy process. Hence, a  L\'evy type process has a natural interpretation as   a  ``process with locally independent increments'', whose characteristic triplet depends on the spatial variable. This justifies the names  \emph{ ``a L\'evy type processes''} and   \emph{``a L\'evy type operator''}.  On the other hand, let a Markov process $X$ be \emph{Feller}, that is,  the respective semigroup $(S_t)_{t\geq 0}$
\begin{equation}\label{Tt}
S_t f(x):= \Ee^x f(X_t)
 \end{equation}
maps the space $C_\infty(\rn)$ of continuous functions vanishing at infinity into itself.
  The  Courr\`ege-Waldenfels  theorem (cf. \cite[Th.~4.5.21]{Ja01}, \cite[Th.~2.21]{BSW}) states that if for a Feller process   the generator of $(S_t)_{t\geq 0}$ is well defined on  $C_c^\infty(\rn)$ (the space of compactly supported infinitely differentiable functions), then on $C_c^\infty(\rn)$ this generator admits  representation (\ref{gener}). Heuristically, this means that a L\'evy type process is a generic form for a Feller processes on $\rn$, and a Feller process on $\rn$ naturally   gives rise to an integral-differential operator  of the form \eqref{gener}.

The converse problem, i.e.  how to show that  an operator of the form \eqref{gener} gives rise to a Feller process and, moreover, to investigate the distribution properties of this process, is  highly non-trivial; this is the topic the current paper is focused on.
Heuristically, the relation between L\'evy type processes and L\'evy processes is similar to that between  diffusions and the Brownian motion. Hence, the problem how to  construct a L\'evy type process given a L\'evy type operator is  similar (but  technically is much more involved)  to the problem of construction   a diffusion process with given coefficients.

There are several ways how to associate a L\'evy type operator  \eqref{gener} to a Markov process. One way is to solve the \emph{martingale problem} for $(L, C_\infty^2(\rn))$, that is, to find a family of probability measures $\Pp^x, x\in \rn$,  such that $\Pp^x(X_0=x)=1$  and  the process
$$
f(X_t)-f(X_0) -\int_0^t Lf(X_s)ds
$$
is a $\Pp^x$-martingale for any  $f\in C^2_\infty(\rn)$, $x\in \rn$.  The difficult part in this problem is to show that the martingale problem is \emph{well posed}, i.e., that the family $\Pp^x, x\in \rn$ is unique. See  \cite{T70}, \cite{T72}, \cite{Km84}, \cite{Ba88}, \cite{Ho94}, \cite{Ho95}, as well as the survey paper \cite{Ba04} and the monograph \cite{Ja05}.
Note that although  the  martingale problem approach is an efficient tool for constructing the process, typically it does not give much information about its intrinsic distribution properties.

Another natural way is based on the general fact that at least in the simplest cases the transition probability density $p_t(x,y)$ of the process $X$ is a  \emph{fundamental solution}  to the Cauchy problem associated with the operator
\begin{equation}\label{A0}
\partial_t - L,
 \end{equation}
see Section \ref{model} below for the definition. In the  parabolic case, i.e. when the integral part in \eqref{gener} is absent, the classical \emph{parametrix method} makes it possible both to construct the fundamental solution $p_t(x,y)$  to (\ref{A0}), and to give the explicit upper and lower bounds for $p_t(x,y)$;  see the monograph by Friedman \cite{Fr64} for  details.  We also refer to the original paper by E.\ Levi \cite{Le07}, and to the paper by W.\ Feller \cite{Fe36}, in which the parametrix construction of the transition probability density  is given for continuous, purely discontinuous, and mixed processes.
Since  for $L$  the \emph{positive maximum principle} holds true, we can conclude  that $p_t(x,y)$ is the transition probability density of a Feller process,  and this process is the unique one associated with  the generator $A=L$. We omit the details, since the same procedure will be discussed in details below in a much more complicated setting.

The goal of our investigation is to extend the approach outlined above to the genuinely ``L\'evy type'' case, where   the diffusion part in \eqref{gener} is absent, but instead the jump part is present: $Q\equiv 0$,  $\mu(x,du)$ is non-trivial. In this case, the structural assumptions on the L\'evy kernel  $\mu(x,du)$ appear to be substantial. The case of $\mu(x,du)$ being comparable, in a sense, to the L\'evy measure of an $\alpha$-stable process $\mu(du)=c|u|^{-\alpha-n}du$ is well studied, see \cite{Dr77},  \cite{DE81}, \cite{Ko89}, \cite{Ko00},  \cite{CZ13}, and an  extensive overview in the monograph \cite{EIK04}. Extending these results to more general classes of L\'evy measures meets new serious difficulties, which we discuss in details below.

Following the line of  the classical ``parabolic''  parametrix method  (cf. \cite{Le07}, \cite{Fr64}), in which the Gaussian kernel is taken as the zero order approximation of the solution to the respective Cauchy problem, the natural idea to develop a ``L\'evy type'' parametrix method is  to take  as the zero order approximation $p_t^0(x,y)$ for the candidate for being the fundamental solution to \eqref{A0} the transition probability density of some L\'evy process. An important feature used in all the aforementioned papers is that the fundamental solution $g_t(x,y)$ to the respective  constant coefficient Cauchy problem satisfies the upper estimate
\begin{equation}\label{fr}
g_t(x,y)\leq C \rho_t^n f(\rho_t (y-x)),
\end{equation}
where  $\rho: (0,1]\to(0,\infty)$ and $f\in L_1(\rn)$ has the meaning of a ``scaling function'' and a ``shape function'', respectively (namely, one has $\rho_t=t^{-1/2}, f(x)=\exp(-\|x\|^2/(2c))$ in the diffusive case and $\rho_t=t^{-1/\alpha}, f(x)=1\wedge \|x\|^{-\alpha-n}$ in the symmetric $\alpha$-stable case).

However, for general L\'evy processes  estimate (\ref{fr}) does not hols, see \cite{KK12a} for a counter-example.
 In \cite{KK12a} and  \cite{K13} it is shown that  a natural analogue of (\ref{fr}) is the so-called  \emph{compound kernel}  upper bound, see Proposition~\ref{aux1} below. An important ingredient in the approach  which we develop in the current paper is   that in the parametrix construction  of the kernel $p_t(x,y)$, \emph{single kernel}-type upper bound (\ref{fr}) can  be successfully replaced by  a  \emph{compound kernel}  upper bound, valid under more mild structural assumptions on the model.

Our construction consists of the following three principal steps. First, we construct the kernel $p_t(x,y)$, which is a \emph{candidate} for being the fundamental solution to the Cauchy problem associated with the operator (\ref{A0}), represented   in the form of convergent series \eqref{sol10}--\eqref{Psi}; see the detailed discussion in Section~\ref{constr}. We  emphasize that in our case the  verification of the fact that the constructed kernel $p_t(x,y)$ is   indeed the fundamental solution hardly could be performed in the classical way described in \cite{Fr64}. The reason for this is that the space derivatives of the zero order approximation $p_t^0(x,y)$ have stronger singularity  at $t=0$, and one cannot prove directly that $p_t(\cdot,y)$ belongs to $C^2_\infty(\rn)$, which is the domain of  $L$.

Our second step is to prove that the constructed kernel $p_t(x,y)$ is the transition probability density of some Markov process, and the restriction of the generator of this process to  $C_\infty^2(\rn)$  equals $L$. The method we use to do this is described in  \cite{KK14a}  (see also \cite{Ku15})\normal\,, and is based on the auxiliary notion of the \emph{approximative fundamental  solution}; see Section~\ref{ver}.

The final step is to identify uniquely the Markov process $X$ obtained before in  terms of the initial operator $L$. While spatial derivatives of $p_t(x,y)$ are hardly controllable, the time derivative is more manageable, which makes it possible to prove that the generator $(A, D(A))$ of $C_\infty(\rn)$-semigroup of $X$ is the closure of $(L, C_\infty^2(\rn))$. This yields   that the martingale problem for $(L, C_\infty^2(\rn))$ is well posed, and shows the uniqueness of the Markov process $X$, constructed in the first two  steps. Also,  we are able to show that $p_t(\cdot,y)$ belongs to the domain of the generator $A$ ($\Leftrightarrow$ the closure of $L$), and that $p_t(x,y)$ is indeed the fundamental solution to the Cauchy problem for
 \begin{equation}\label{A}
\partial_t - A,
 \end{equation}
which justifies our parametrix construction. Thus,  starting with the Cauchy problem for (\ref{A0}), we construct the fundamental solution for (\ref{A}), where $A$ is the closure of $L$.

As an application of the estimates constructed for the kernel,   we give  the  necessary and (separately) sufficient conditions  for a  finite Borel measure to  belong to the Kato  and Dynkin classes with respect to $p_t(x,y)$.

Let us give a brief overview of other   existing results.

 In the L\'evy case  the transition probability density is just the inverse Fourier transform of the charactristic function. This allows a  lot of possibilities to estimate this transition probability density, see, for example,    \cite{PT69}, \cite{H94}, \cite{H03}, \cite{W07} for the asymptotic behaviour of an $\alpha$-stable transition probability density,  and \cite{KR15},   \cite{KSt13a}--\cite{KSt13c}\normal\,,   \cite{RS10},    \cite{Mi15},\normal\, \cite{BGR14}, \cite{KK12a}, \cite{K13},   \cite{St10a}--\cite{St16}, \cite{CGT15}\normal\,  for the L\'evy case. Of course, this list of publications is far from being complete.

In  \cite{BJ07}  the case of the  fractional Laplacian perturbed by a gradient is treated, see also  \cite{BJS12}, \cite{BS12} for the  kernel estimates.   The verification procedure presented in \cite{BJ07}  shows that the integro-differential operator is the weak generator of the respective semigroup.  In \cite{KS14} the case of a singular perturbation of the fractional Laplacian is considered, and a different approach is used for the verification; see also \cite{CW13},    for an approach which relies on \cite{BJ07} and the martingale problem,  as well as \cite{CZ13} and the references therein.\normal\,
  In  \cite{Po94} and \cite{PPo95} the authors constructed  the transition probability density of the process  which is the weak solution to the
SDE driven by a symmetric   $\alpha$-stable process with a drift.
We refer to \cite{FP10}, \cite{DF13}, and \cite{KK14a}, in which the gradient perturbations of an $\alpha$-stable like operator with $0<\alpha<1$ are investigated.

Another approach to study the fundamental solution to the respective Cauchy problem
    relies on  a different  version of the parametrix
method and the symbolic calculus technique, which allows to prove  the existence of
the fundamental solution, and to construct it in the form of converging in a certain
sense series. This approach uses the Hilbert space methods, and is developed in \cite{Ts74}, \cite{Iw77}, \cite{Ku81},
\cite{Ho98a}, \cite{Ho98b}, \cite{Ja02}, 
\cite{B\"o05} and  \cite{B\"o08}.

 There is a large group of results devoted to the estimation of the transition probability density of a Markov process,  associated with a Dirichlet form of a certain type. Under the assumption that  the
 jump intensity measure of a Markov process is absolutely continuous and possesses certain regularity properties, estimates on the
  transition probability density   are obtained in \cite{CKS87}, \cite{CK08},
\cite{CKK08}, \cite{CKK10}, \cite{BBCK09}, \cite{BGK09}, \cite{Mi12}; of course,
this list is far from being complete. The approach used in the above
papers  relies on the  Dirichlet form technique and the Harnack principle. Note that in these papers  the initially given object is a regular Dirichlet form, which  already assumes  the \emph{existence} of the related Markov process.

The  paper is organized as follows. In  Section~\ref{set}
we set the notation,   outline the method,  and formulate the results.  The construction of the parametrix series is performed in Sections~\ref{aux-res}--\ref{pr-main}. Proofs of the continuity and smoothness properties are given in Section~\ref{pr-cont}. Section~\ref{ver} is devoted to the verification  procedure. The uniqueness is studied in Section~\ref{pr-gen}.  Diagonal  and lower bounds for the constructed  fundamental solution are given in Section~\ref{up-lo}. Finally, Section~\ref{appl} is devoted to the proof of the application result, that is, using the structure of the upper and lower bounds on $p_t(x,y)$, we provide the necessary and (separately) sufficient conditions for a measure to be in the Kato and Dynkin classes.

\section{The main results: outline and formulation} \label{set}

\subsubsection*{Notation}

 For functions $f$, $g$ we mean by $f\asymp g$
that there exist some constants $c_1, c_2>0$ such that $c_1 f(x)
\leq g(x) \leq c_2 f(x)$ for  all $x\in \rn$. By $x\cdot y$ and
$\|x\|$  we denote, respectively, the scalar product and the norm in
$\rn$;  $\Sss^n$ denotes a unit sphere in $\rn$. We denote by $c_i$,
$c$, $C$, etc., arbitrary positive constants. Denote by $(f*g)(t,x,y)$, $(f\star g)(t,x,y)$, $(F*G)_t(du)$, $(F\star G)_t(du)$ the respective convolutions of functions $f(t,x,y)$, $g(t,x,y)$, and of kernels  $F_t(du)$ and $G_t(dv)$:
 $$
(f* g)(t,x,y):= \int_\rn f(t,x,z)g(t,z,y)dz, \quad (f\star g) (t,x,y):= \int_0^t \int_\rn f(t-s,x,z)g(s, z,y)dzds,
 $$
$$
(F* G)_t (du)= \int_\rn F_t(du-z) G_t(dz), \quad (F\star G)_t(du)=\int_0^t \int_\rn
F_{t-s}(du-z) G_s(dz)ds.
$$
By $B_b(\rn)$, $C_\infty^k (\rn)$,  we denote, respectively, the sets of bounded Borel functions,  and the set of $k$-times differentiable functions, vanishing at infinity together with their derivatives. By $\|\cdot\|_\infty$ we  define the $\sup$-norm in $C_\infty$.

\subsection{The model and the outline of the method}\label{model}

  In this section we describe in detail three steps of our approach, indicated in the Introduction.

 Let $L$ be an operator of the form \eqref{gener} but with $Q\equiv 0$. Below we specify the assumptions on the drift and the kernel. \normal
In the first part we construct a  candidate for the  fundamental solution to the Cauchy problem
\begin{equation}\label{L_ful}
\partial_t-L,
\end{equation}
i.e. such a  function $p_t(x,y)$ that
\begin{equation}\label{L_delta}
p_t(x, \cdot)\to \delta_x\quad \text{ as $t\to 0+$, $x\in \rn$,}
\end{equation} and
\begin{equation}\label{L_fund}
\Big(\prt_t-L\Big)p_t(x,y)=0, \quad t>0, \quad x,y\in \rn.
\end{equation}
In order to simplify the further exposition, let us briefly outline the parametrix  construction, see   \cite[p.310--311]{Fr64},  or  \cite[p.144--145]{Ja02} for more information.

Consider \emph{some}
 approximation $p_t^0(x,y)$ of $p_t(x,y)$, and denote by $r_t(x,y)$ the residue with respect to this approximation, that is,
\begin{equation}\label{sol10}
p_t(x,y)=p_t^0(x,y)+ r_t(x,y).
\end{equation}
Define
\begin{equation}\label{Phi}
\Phi_t(x,y):=\Big(L-\prt_t\Big)p_t^0(x,y),\quad t>0, \quad x,y\in\rn.
\end{equation}
Observe that since  $p_t(x,y)$ is aimed  to be the fundamental solution to the Cauchy problem for the operator  (\ref{L_ful}), we should have
$$
\Big(\prt_t-L\Big)r_t(x,y)=\Phi_t(x,y).
$$
Therefore,
$$
 r_t(x,y)=(p\star \Phi)_t(x,y),
$$
 which by \eqref{sol10} allows us  to write the equation for  $r_t(x,y)$:
$$
r_t(x,y)=(p^0\star \Phi)_t(x,y)+(r\star \Phi)_t(x,y).
$$
The formal solution to this equation is given by the convolution
\begin{equation}\label{r}
r=p^0\star \Psi,
\end{equation}
where
\begin{equation}\label{Psi}
\Psi_t(x,y)= \sum_{k=1}^\infty \Phi^{\star k}_t(x,y).
\end{equation}
We can chose the zero order approximation $p_t^0(x,y)$ in the following way.

Consider the operator
 \begin{equation}\label{lz}
 L^zf(x):=a(z) \cdot  \nabla f(x)+ \int_\rn \big(f(x+u)-f(x)-u \cdot \nabla f(x) \I_{\{\|u\|\leq 1\}}\big)\mu(z,du),
 \end{equation}
 where $f\in C_\infty^2(\rn)$.  It is known that $(L^z, C_\infty^2(\rn))$ extends to the generator of  a semigroup corresponding to a L\'evy process, which under condition \textbf{A1} (see below) possesses the transition probability density $\mathfrak{p}_t^z(x)$.  Note that $\mathfrak{p}_t^z(y-x)$ is the fundamental solution to a Cauchy problem  for the operator
$$
\partial_t-L^z,
$$
  see \cite[Ex.~2.7.14]{Ja02}.   Put
\begin{equation}\label{pto}
p_t^0(x,y):= \mathfrak{p}^z_t(y-x)\big|_{z=y}.
\end{equation}
We prove that under  such a  choice of the zero order  approximation $p_t^0(x,y)$, the series in \eqref{Psi} indeed converge, and that the expression \eqref{sol10} is well defined.

  On the second step we associate with the constructed kernel  $p_t(x,y)$ a Markov process. The keystone in this step is the usage of an auxiliary object, which we call the \emph{approximate fundamental solution}, which is  a certain approximation $p_{t,\epsilon}(x,y)$ of the constructed kernel $p_t(x,y)$.
Using the expression for $p_t(x,y)$ and the estimates on $\Phi_t(x,y)$ and $\Psi_t(x,y)$, obtained in the first step, we show that for the operators $S_{t,\eps}$ with kernel $p_{t,\eps}(x,y)$  the following statements hold true:
 \begin{itemize}
  \item[(a)]  For $f\in C_\infty(\rn)$,
  $$
  \lim_{\eps\to 0}\|S_{t,\eps}f- S_t f\|_\infty=0,
  $$
  uniformly on compact subsets of $(0,\infty)$;

  \item[(b)] For $S_{t,\eps} f$, $f\in C_\infty(\rn)$,  $\eps>0$, identity (\ref{L_fund}) turns into the approximative identity
      $       \Big(\partial_t -L\Big)S_{t,\eps} f$,
      and
       $$
  \lim_{t,\eps\to 0}\|S_{t,\eps}f- S_t f\|_\infty=0;
  $$
\end{itemize}
here
\begin{equation}\label{t-semi}
S_tf(x):= \int_\rn f(y) p_t(x,y)dy,\quad t>0, \quad x\in \rn.
\end{equation}
These properties of $S_{t,\eps}$ allow us to develop a version of the positive  maximum principle (see  \cite[p.165]{EK86} or \cite[Cor.~4.5.14]{Ja01}  for the classical positive maximum principle), which in turn enables us  to show that $p_t(x,y)$ is the transition probability density of a Markov process, which is a solution to the martingale problem for $(L, C_\infty^2(\rn)$. In particular, the family of operators $(S_t)_{t\geq 0}$ forms a semigroup,  related to $X$ by
$$
S_tf(x)= \Ee^x f(X_t), \quad f\in C_\infty^2(\rn).
$$

The third step is  devoted to the uniqueness problem  for the constructed process.
We show that the generator $A$ of the semigroup $(S_t)_{t\geq 0}$ coincides on $C_\infty^2(\rn)$ with $L$.  Further,  we employ the properties of the derivative $\prt_t S_{t,\eps} f$:
\begin{itemize}
\item[(c)] For any $f\in C_\infty(\rn)$,
 $$
 \lim_{\eps\to 0}\|\partial_t S_{t,\epsilon}f- \prt_t S_t f\|_\infty=0,
 $$
   uniformly on compact subsets of $(0,\infty)$.
\end{itemize}
This property together with (b) allows to control $L$ on $S_{t,\eps} f$,  $f\in D(A)$  and  show, that $A$ is the closure of $L$ in $C_\infty$. Consequently, the process $X$ constructed in the previous step is the unique solution to the above martingale problem.

Finally, the property similar to (c)  holds also for the kernels $p_{t,\eps}(x,y)$ and $p_t(x,y)$:

 \begin{itemize}
 \item[(d)]   $\prt_t p_{t,\eps} (x,y)$ approximates $\prt_t p_t(x,y)$ as $\eps\to 0$, uniformly on compact subsets of $(0,\infty)\times \rn\times \rn$.
 \end{itemize}
 Therefore, using similar  machinery we show that  $p_t(\cdot,y)\in D(A)$  for any fixed $y$, and
that  $p_t(x,y)$  is the fundamental solution to the Cauchy problem for $\partial_t -A$.

\subsection{Main results}\label{main}

Consider
\begin{equation}\label{qu110}
q(\xi) :=\int_\rn (1-\cos (\xi \cdot u))\mu(du),
\end{equation}
where $\mu(du)$ is a  L\'evy measure, i.e. a Borel measure satisfying $\int_\rn (\|u\|^2 \wedge 1 )\mu(du)<\infty$, and define
    \begin{equation}
    q^U(\xi):=\int_\rn [(\xi \cdot u)^2\wedge 1]\mu(du),\quad q^L(\xi):=\int_{|u\cdot \xi|\leq 1} (\xi \cdot u)^2\mu(du).\label{psipm}
    \end{equation}
    The function $q(\xi)$  possesses the L\'evy-Khinchin representation, and thus is the characteristic exponent of a L\'evy process.
It can be shown  (cf. \cite{KK12a})   that the functions $q^L(\xi) $ and $q^U(\xi)$ satisfy
 \begin{equation}
(1-\cos 1) q^L(\xi)\leq q(\xi) \leq 2q^U(\xi). \label{psipm1}
\end{equation}
 Note that in \eqref{qu110} and \eqref{psipm} we do not assume  $\mu$ to be  symmetric.  \normal Suppose  that the measure $\mu$ satisfies the regularity assumption given below.

\begin{itemize}
\item[\textbf{A1.}]   There exists $\beta>1$ such that
$$
\sup_{l\in \Sss^n}q^U(r l)\leq \beta\inf_{l\in \Sss^n} q^L(r l)\quad \text{ for all $r>0$  large enough.}
$$
\end{itemize}

In what follows, we denote  \begin{equation}\label{al}
  \alpha:= 2/\beta.
  \end{equation}
This notation is motivated by the  particularly important example of  a symmetric $\alpha$-stable L\'evy measure $\mu(du):=c(\alpha) \|u\|^{-n-\alpha} du$, $\alpha\in (0,2)$: direct calculations show that in this case \textbf{A1} holds true with $\beta=2/\alpha$. Note also that  for any L\'evy measure $\mu$ satisfying \textbf{A1} the respective  L\'evy exponent  $q$ admits a polynomial lower bound (see (\ref{qal}) below), which for the  symmetric $\alpha$-stable L\'evy measure becomes  an identity.

Through the paper we assume that  the kernel $\mu(x,du)$ is of the form
\begin{equation}
\mu(x,du) = m(x,u)\mu(du), \label{mx}
\end{equation}
where $m(x,u)$ is some  positive  measurable function.   We assume that the function $m(x,u)$ and the drift coefficient  $a(x)$ satisfy the assumptions below.

\begin{itemize}
\item[\textbf{A2.}]  The functions  $m(x,u)$ and $a(x)$ are measurable, and satisfy  with some constants $b_1,\, b_2,\, b_3 >0$ the inequalities
    $$
    b_1\leq m(x,u)\leq b_2,\quad |a(x)|\leq b_3, \quad  \quad x,u\in \rn.
    $$

\item[\textbf{A3.}]  There exist   constants $\gamma \in (0,1]$ and $b_4>0$ such that
\begin{align}\label{M2b}
 |m(x,u)-m(y,u)| + \|a(x)-a(y)\| \leq  b_4( \|x-y\|^\gamma \wedge 1),\quad u,\, x,\,y \in \rn.
\end{align}

\item[\textbf{A4.}]  In  the case $\alpha\in (0,1]$ we  assume that $a(x)=0$ and  the kernel $\mu(x,du)$ is symmetric with respect to $u$ for all $x\in \rn$.
\end{itemize}

  Below we state the first main result of our paper.

\begin{theorem}\label{t-main1}
 Suppose that assumptions \textbf{A1} -- \textbf{A4} are satisfied, and the function $p_t^0(x,y)$ is given by \eqref{pto}. Then

\begin{itemize}
\item[a)] The  function $p_t(x,y)$  introduced in  \eqref{sol10} -- \eqref{Psi}  is well defined in the sense that the series \eqref{Psi} and  \eqref{r} converges absolutely for any $t>0$, $x,y\in \rn$, uniformly on compact subsets of $(0,\infty)\times \rn\times \rn$;

\item[b)]  The function $p_t(x,y)$ is continuous on $(0,\infty)\times \rn\times \rn$.
\end{itemize}

\end{theorem}
Next we associate the constructed  function $p_t(x,y)$ with the initial operator $L$.  To make the structure the most transparent, we do this in two steps:  we prove that $p_t(x,y)$ a  transition probability density of \emph{some} Markov process, and then  show that the $C_\infty(\rn)$-generator of this process is an extension of $(L, C^2_\infty(\rn))$. The second statement means that  the semigroup (\ref{t-semi}) with $p_t(x,y)$ defined by (\ref{sol10}) is in fact \emph{unique} Feller semigroup  associated with the operator $L$.

\begin{theorem}\label{t2}  The family of operators \eqref{t-semi}
forms  a strongly continuous conservative semigroup of non-negative operators  on $C_\infty(\rn)$, which in turn defines a (strong) Feller Markov process $X$. Further, the set  $C^2_\infty(\rn)$ belong to the domain $D(A)$ of the generator $A$ of this semigroup, and
 $$
 Af(x)=Lf(x) \quad\text{for} \quad  f\in C_\infty^2(\rn),
 $$
 that is, $(A,D(A))$ is an extension of  $(L, C^2_\infty(\rn))$.
\end{theorem}

\begin{theorem}\label{t3}
\begin{itemize}
\item[a)] The generator $(A,D(A))$ is the closure of $(L, C^2_\infty(\rn))$.

\item[b)] The function $p_t(\cdot ,y)$ belongs to the domain $D(A)$ of $A$, and  is the fundamental solution to the Cauchy problem  for the operator $\partial_t - A$.
    \end{itemize}
\end{theorem}
The first statement of Theorem~\ref{t3} allows us to show the uniqueness of the solution to the martingale problem for $(L, C_\infty^2(\rn))$.
\begin{theorem}\label{mart}
The Markov process $X$ constructed in Theorem~\ref{t3} is the unique solution  to the martingale problem for $(L, C_\infty^2(\rn))$.
\end{theorem}
Finally,  we give the upper and lower estimates for the constructed function $p_t(x,y)$ and its time derivative.

Let
\begin{equation}
q^* (r) :=\sup_{l\in \Sss^n}q^U(r l), \quad   r>0,   \label{qu}
\end{equation}
  and define
\begin{equation}
\rho_t:= \inf\{ r:\quad  q^*(r) =1/t\}, \quad   t\in (0,T].  \label{rho1}
\end{equation}
   Since the function $q^*(r)$ is continuous and $\lim_{r\to \infty} q^*(r)=\infty$, the function $\rho_t$, $t\in (0,T]$, is well defined for any $T>0$. \normal

In  \cite{K13}, see also \cite{KK12a}, we show  that condition \textbf{A1} implies   for  $r$ large enough the lower estimate
\begin{equation}
q^*(r)\geq c_0r^\alpha,
\label{qal}
\end{equation}
which in turn  implies for any $T>0$ the upper bound
\begin{equation}\label{growth2}
\rho_t \leq c_1 t^{-1/\alpha}, \quad t\in (0,T].
\end{equation}
Note that    for any $c>1$ we have $q^U (c \xi)\leq (c^2\wedge 1) q^U(\xi)$,   implying  $q^*(r)\leq c_2 r^2$, $r\geq 1$,  then $\rho_t\geq c_3 t^{-1/2}$, $t\in (0,T]$. Denote by $\sigma\in [\alpha,2]$ the minimal value for which there exists $c_4>0$ such that
\begin{equation} \label{growth3}
\rho_t\geq c_4 t^{-1/\sigma}, \quad t\in (0,T].
\end{equation}

\medskip

Denote by $f_{up}$ and $f_{low}$ the functions  of the form
\begin{equation}
f_{up}(x):=d_1 e^{-d_2 \|x\|},\quad  f_{low}(x):= d_3 (1-d_4\|x\|)_+, \quad    x\in \rn, \label{flu}
\end{equation}
where  $d_i>0$, $1\leq i\leq  4$, are some constants which are yet to be chosen.

\begin{theorem}\label{t4}
 For any $T>0$ \normal there exist  constants $d_i>0$, $1\leq i\leq 4$,  and a family of   sub-probability measures \,  $\{Q_t, \, t\geq 0\}$,  such that   $p_t(x,y)$ satisfies the upper and lower estimates
\begin{equation}\label{eqII}
\rho_t^nf_{low}(\rho_t(y-x)) \leq p_t(x,y)\leq  \rho_t^n \big(f_{up}(\rho_t \cdot) * Q_t\big)(y-x), \quad t\in (0,T], \quad x,y\in \rn.
\end{equation}
 where $f_{low}$ and $f_{up}$ are the functions of the form \eqref{flu} with  constants $d_i$, $1\leq i\leq 4$.\normal
\end{theorem}
\begin{theorem}\label{t5}
\begin{itemize}
\item[1.] There exists  $\partial_t p_t(x,y)$, which is continuous in   $(t,x,y)\in (0,\infty)\times \rn\times \rn$.

\item[2.]   For any $T>0$  \normal there exist constants $\tilde d_1, \tilde d_2>0$ and   a family of   sub-probability measures \,    $\{\tilde{Q}_t, \,t\geq 0\}$,  such that
   $$
   |\prt_tp_t(x,y)|\leq t^{-1}\rho_t^n \big(f_{up}(\rho_t \cdot) * \tilde{Q}_t\big)(y-x),\quad  t\in (0,T], \quad x,y\in \rn.
   $$
 where $f_{up}$  is of the form  \eqref{flu} with constants $\tilde d_1,\tilde d_2$.\normal
\end{itemize}

\end{theorem}

To demonstrate an application of  the above results, we need a bit more preparations.

Recall that a functional $\varphi_t$ of a strong Markov process $X$ is called a $W$- functional, if it is additive, positive, continuous, almost surely homogeneous, and
$$
v_t(x):= \Ee^x \varphi_t<\infty;
$$
 in this case the function $v_t(x)$ is called the characteristic of $\varphi_t$, see \cite[$\S6.11$]{Dy65}. By \cite[Th.~6.3]{Dy65}, the characteristic determines the $W$- functional uniquely up to equivalence. On the other hand,  \cite[Th.~6.3]{Dy65} gives a way how to check that a given function is a characteristic of  some $W$-functional.

Recall (cf. \cite{KT07},  \cite{AM92}) that
  a  Borel measure $\varpi$ is said to belong to

  \begin{itemize}
  \item[i)]
the Kato class $S_K$ with respect to  $p_t(x,y)$,  if
  \begin{equation}
    \lim_{t\to 0} \sup_{x\in \rn}  \int_0^t \int_\rn p_s(x,y)\varpi(dy)ds=0;  \label{K1}
    \end{equation}
\item[ii)] the  Dynkin class $S_D$ with respect to  $p_t(x,y)$,  if  there exists  $t>0$ such that
    \begin{equation}
    \sup_{x\in \rn}  \int_0^t \int_\rn p_s(x,y)\varpi(dy)ds<\infty. \label{W1}
    \end{equation}
    \end{itemize}
 Clearly, $S_K\subset S_D$. By  \cite[Th.~6.6]{Dy65},   condition $\varpi\in S_K$    implies  that the function
 \begin{equation}
\chi_t(x):= \int_0^t \int_\rn p_s(x,y)\varpi(dy)ds\label{chit}
\end{equation}
is the characteristic of some $W$-functional $\varphi_t$,
provided that the mapping $x\mapsto \chi_t(x)$ is measurable for each $t\geq 0$.
Thus, to prove that  $\chi_t(x)$ is the characteristic of some $W$- functional of $X$,    we need to check whether the measure $\varpi$ from \eqref{chit} belongs to the Kato class with respect to $p_t(x,y)$. As an accompanying result, we get the condition  under which $\varpi$ belongs to the respective Dynkin class.

\begin{remark}
Up to our knowledge there are not many results on the necessary and sufficient conditions when a measure is in the Kato class.   In the case of a symmetric $\alpha$-stable process, $\alpha\in (0,2)$, and a relativistic $1/2$-stable process, 
these conditions  are stated   in  \cite{Z91}, see also \cite{FOT94}.
In the case of $n$-dimensional Brownian motion there is a one-to-one correspondence between the class of $W$-functions and so-called $W$-measures,  see  \cite[Th.~8.4]{Dy65}; in our notation this theorem means that every measure from the Dynkin class is in one-to-one correspondence with a  $W$-functional.   An example of a measure which for a Brownian motion belongs to the class $S_D$ but not to the  $S_K$, can be found e.g. in    \cite{Ku09}.
\end{remark}

In the theorem below we present the necessary and (separately) sufficient conditions when a measure belongs to the Kato and Dynkin classes with respect to $p_t(x,y)$.
 \begin{theorem}\label{loc-time}
 Let $\varpi$ be a finite Borel measure on $\rn$.

\begin{itemize}

\item[a)]  For $\varpi\in S_D$ with respect to  $p_t(x,y)$ it is sufficient that
\begin{equation}\label{d-suf}
 \int_0^\delta \frac{\sup_{x\in \rn} \varpi\{ y:\, \|x-y\|\leq s\} }{s^{n+1} q^*(1/s)}ds<\infty, \quad \text{for some $\delta>0$},
\end{equation}
and necessary, that
\begin{equation}\label{d-nes}
\sup_{x\in \rn} \int_0^\delta \frac{\varpi\{ y:\, \|x-y\|\leq s\} }{s^{n+1} q^*(1/s)}ds<\infty, \quad \text{for some $\delta>0$}.
\end{equation}

\item[b)]  For $\varpi\in S_K$ with respect to  $p_t(x,y)$ it is sufficient that \eqref{d-suf} holds true,
and necessary, that
\begin{equation}\label{k-nes}
\lim_{\delta\to 0} \sup_{x\in \rn}\int_0^\delta \frac{\varpi\{ y:\, \|x-y\|\leq s\} }{s^{n+1} q^*(1/s)}ds=0.
\end{equation}
\end{itemize}
\end{theorem}

\section{Construction of the parametrix series. Proof of Theorem~\ref{t-main1}}\label{constr}

\subsection{Well-definiteness of $p_t^0(x,y)$ }\label{aux-res}

It is known that for any fixed $z\in \rn$  the operator $(L^z, C_\infty^2(\rn))$ (see  \eqref{lz}) extends to the $C_\infty$- generator of a Feller semigroup which corresponds to the L\'evy process $X_t^z$ with  characteristic function
$$
\Ee e^{i\xi \cdot X_t^z} = e^{-t q(z,\xi)},
$$
where
\begin{equation}
q(z,\xi) := -i a(z) \cdot \xi +\int_\rn (1-e^{i\xi \cdot u}+i\xi \cdot u \I_{\{\|u\|\leq 1\}})\mu(z,du). \label{qu11}
\end{equation}
Note that due to  condition \textbf{A2} the kernels $\{\mu(z,du),\,\,z\in\rn\}$ are comparable in the sense that for any $z, y\in \rn$ and any Borel subset $A\subset\rn\backslash \{0\}$ we have $\mu(z,A)\asymp \mu(y,A)$, implying that
\begin{equation}
\RR q(z,\xi)\asymp \RR q(y,\xi) \quad \quad \text{ for all $z,y,\xi\in \rn$}.
 \end{equation}
 Condition \textbf{A3}  implies that $\RR q(z,\xi)$ is continuous in $z$.
 Condition   \textbf{A1} together with \textbf{A2} implies  (cf. \cite{KK12a}, \cite{K13})  that
\begin{equation}
\min_z \RR q(z,\xi)\geq c\|\xi\|^{\alpha}\quad \text{ for large $\|\xi\|$}, \label{growth}
\end{equation}
where $\alpha=2/\beta$. Thus,   for any fixed $z$ the  process $X_t^z$   admits a transition probability density, which we denote    by
 $\mathfrak{p}^z_t(x)$.  Note that by \eqref{growth} we
 can write $\mathfrak{p}^z_t(x)$ as
\begin{equation}\label{qtz}
\mathfrak{p}^z_t(x)= (2\pi)^{-n} \int_\rn e^{-ix\cdot \xi -tq(z,\xi)}d\xi.
\end{equation}
Thus, the function $p_t^0(x,y)$ given by \eqref{pto} is well defined, and  $p_t^0(\cdot,y)\in C_b^\infty(\rn)$.

\subsection{Estimate for $\Phi_t(x,y)$}\label{Ph1}

 In this subsection we derive the upper bound for $\Phi_t(x,y)$, see Lemma~\ref{Phi-up}. In order to do this, we introduce some notation and state the auxiliary propositions, the proofs of which we defer to Appendix A.

Let
\begin{equation}
\Lambda_t(du):=t  \mu(du)\I_{\{\rho_t\|u\|>1\}}.    \label{lam}
\end{equation}
\begin{proposition}\label{lem1}
For  any $T>0$  we have $\Lambda_t(\rn) \leq   n^2$, $t\in [0,T]$.
\end{proposition}

\begin{proposition}\label{cor11} For any $\lambda\in [0,\alpha)$, $T>0$, we have
\begin{equation}
 \rho_t^\lambda \int_\rn (\|u\|^\lambda\wedge 1) \Lambda_t(du)  \leq C, \quad t\in [0,T].
 \end{equation}
\end{proposition}
\normal
Define the probability measure
\begin{equation}
   P_t(du) :=e^{-\Lambda_t(\rn)} \sum_{k=0}^\infty \frac{1}{k!} \Lambda_t^{*k}(du). \label{Poist}
    \end{equation}
Let $\alpha$ be the parameter defined in \eqref{al}, and  let  $\gamma\in (0,1]$ be  the parameter of H\"older continuity from \textbf{A3}.  Fix some $\epsilon\in (0,\alpha)$, and put
\begin{equation}\label{kap}
\kappa:=
\begin{cases}
\gamma, &   \text{if}  \quad  \gamma\in (0,\alpha)\\
\alpha-\epsilon, &   \text{if} \quad \gamma\geq \alpha.
\end{cases}
\end{equation}
  Note that by definition $\kappa>0$.
Put
 \begin{equation}\label{ptg}
P_{t,\kappa}(du):=
\big(1+\rho_t^\kappa \big( \|u\|^\kappa \wedge 1\big)\big)  P_t(du).
\end{equation}

\begin{proposition}\label{ptkappa}
For any $T>0$ there exists $C\in (0,\infty)$ such that $P_{t,\kappa}(\rn)\leq C$, $t\in [0,T]$.
\end{proposition}

 Finally, define
\begin{equation}
G_{t}(du):=c_0\Big( P_{t,\kappa}(du)+   P_t \normal * P_{t,\kappa}(du)\Big),\label{G1}
\end{equation}
where the constant $c_0>0$ is chosen in such a way that   $G_t(\rn)\leq 1$ for all $t\in [0,T]$. Such a choice of $c_0$ is possible due  to   Proposition~\ref{ptkappa}.\normal

Now we  are ready to state the main result of this section. Recall that the parameter $\sigma$ was defined in \eqref{growth3}.
\begin{lemma}\label{Phi-up}
For any $T>0$ there exist constants $C, b>0$  such that
\begin{equation} \label{ld11}
\big|   \Phi_t(x,y)\big|\leq  C t^{-1+\eta} \big(g_t *
G_{t}\big)(y-x),\quad  t\in (0,T], \, x,y\in \rn, \normal
\end{equation}
where   $\eta=\frac{\kappa}{\sigma}\wedge \big( 1+ \frac{\kappa-1}{\alpha}\big)$,\normal
\begin{equation}\label{gt10}
g_t(x):= \rho_t^n  e^{- b \rho_t\|x\|}
\end{equation}
  with some $b>0$\normal, and  $\{G_{t} (\cdot),\, t>0\}$  is the  family of   sub-probability    measures  given by \eqref{G1}.
\end{lemma}
In the proof of this lemma we use several  auxiliary statements which we formulate below.

Take $f_{up}$ of the form \eqref{flu}, and put
  \begin{equation}\label{ft}
f_t(x):= \rho_t^n \big( f_{up}(\rho_t \cdot )*  P_t\big)(x). 
\end{equation}

\begin{proposition}\label{aux1}
Suppose that  conditions \textbf{A1},  \textbf{A2}  and \textbf{A4} hold true.

Then
  for any $k=k_1+\ldots+k_n\geq 0$, $T>0$,  there exist constants $A_k,a_k>0$,  such that
 \begin{equation}\label{der-est0}
\Big| \frac{\partial^k}{\partial x_1^{k_1} \ldots \partial x_n^{k_n}} \mathfrak{p}^z_t(x)\Big|\leq   \rho_t^k f_t(x), \quad t\in (0,T],\quad  x,z\in \rn,
\end{equation}
where $f_t(x)$ is the function of the form \eqref{ft}, and the function $f_{up}$ in the definition of $f_t$ is of the form \eqref{flu}, with constants $A_k$, $a_k$
in the place of $d_1$, $d_2$, respectively.

 In particular,
 \begin{equation}\label{der-est}
\Big| \frac{\partial^k}{\partial x_1^{k_1} \ldots \partial x_n^{k_n}} p_t^0(x,y)\Big|\leq   \rho_t^k f_t(y-x), \quad t\in (0,T],\quad  x,y\in \rn.
\end{equation}
\end{proposition}
\begin{proposition}\label{aux2} Suppose that  conditions \textbf{A1},  \textbf{A2}  and \textbf{A4} hold true.
 For $T>0$ there exist $d_3,d_4>0$ such that
$$
\mathfrak{p}^z_t(x) \geq \rho_t^n  f_{low}(\|x\|\rho_t), \quad x,z\in \rn, \, t\in (0,T].
$$
where $f_{low}$ is of the form \eqref{flu} with these constants $d_3$ and $d_4$. In particular, $$
p_t^0(x,y)\geq \rho_t^n  f_{low}(\|  y-x\normal\|\rho_t), \quad x\in \rn, \, t\in (0,T].
$$
 \end{proposition}
\normal
\begin{remark}\label{ak}
Proceeding as in  \cite{K13} and \cite{KK12a} one can show that $0<a_k\leq a_{k-1}$, $k\geq 1$.
\end{remark}
The proof of this  proposition repeats line by line the proof of a similar statement in \cite{K13}, see also \cite{KK12a}. The only difference is that we need to check, using conditions \textbf{A2}--\textbf{A4},  that the required estimates obtained in \cite{K13} hold true uniformly in $z$.   We omit the details.

\begin{proposition}\label{lem1-1}
For any  $\theta\in (0,1)$, $T>0$,    we have
\begin{equation}\label{ft2}
  (\|x\|^\kappa \wedge 1)   f_t(x) \leq C \rho_t^{-\kappa}\big(  g_{t,\theta}    * P_{t,\kappa}\big)(x),\quad t\in (0,T], \quad x\in \rn,
\end{equation}
where  $\kappa$ is defined in \eqref{kap},
 $P_{t,\kappa}(dw)$ is defined in \eqref{ptg},   $f_t$ is of the form \eqref{ft} with some $f_{up}$, and
 \begin{equation}\label{gt}
 g_{t,\theta}(x)= \rho_t^n f_{up}(\theta \rho_t x).
 \end{equation}
 \end{proposition}

\begin{proof}[Proof of Lemma~\ref{Phi-up}]
By the definition of $p_t^0(x,y)$, for any $y\in\rn$ we have
$$
[\prt_t -L^y_x(D)] p^0_t(x,y)=0.
$$
Then
    \begin{equation}
    \begin{split}
       \Phi_t(x,y)&=[L(x,D)-L^y_x(D)]p_t^0(x,y)
    \\&
=\big(a(x)-a(y)\big) \cdot \nabla  p_t^0(x,y)
\\&
\quad +\int_\rn [p_t^0(x+u,y)-p_t^0(x,y) -u\cdot  \nabla  p_t^0(x,y)\I_{\{\|u\|\leq 1\}}][m(x,u)-m(y,u)] \mu(du)
    \\&
    =\big(a(x)-a(y)\big) \cdot \nabla p_t^0(x,y)
    \\&
\quad +\left[ \int_{\rho_t \| u\|\leq 1}+\int_{\rho_t\| u\|>1} \right] [p_t^0(x+u,y)-p_t^0(x,y) -u\cdot  \nabla  p_t^0(x,y)\I_{\{\|u\|\leq 1\}}]\\
&\quad \quad \cdot [m(x,u)-m(y,u)] \mu(du)
    \\&
    =:J_1+J_2+J_3.\label{l10}
    \end{split}
    \end{equation}
We estimate the terms  $J_i$, $i=1,2,3$, separately.   In what follows,  $f_t$ and    $g_{t,\theta}$   are the functions appearing in  Propositions~\ref{aux1} and \ref{lem1-1}.

\medskip

Note that by \textbf{A4} we have $J_1=0$ for $\alpha\in (0,1]$. For $\alpha\in (1,2)$ we have by  \eqref{der-est} and \textbf{A3} the estimates
 \begin{equation}\label{lj1}
\begin{split}
|J_1|\leq  \sqrt{n} \normal\|a(x)-a(y)\| \rho_t f_t(y-x)
&\leq c_1  (\|y-x\|^\gamma \wedge 1) \rho_t f_t(y-x).
\end{split}
\end{equation}
Using Proposition~\ref{lem1-1}   we obtain
\begin{align*}
(\|y-x\|^\gamma \wedge 1) f_t(y-x)  \leq    (\|y-x\|^\kappa \wedge 1)  f_t(y-x)
\normal\leq c_2 \rho_t^{-\kappa} \big(g_{t,\theta} * P_{t,\kappa}\big)(y-x),
\end{align*}
where $\kappa$  is defined  in \eqref{kap}, the semigroup  $P_{t,\kappa}(dw)$ is defined in \eqref{ptg}, and   $\theta\in (0,1)$ is some constant.
\,  Note that since $\alpha\in (1,2)$, we have  $\kappa=\gamma\leq 1$. \normal
Using \eqref{growth2}, we derive
$$
\rho_t^{1-\kappa}\leq c t^{-\frac{1-\kappa}{\alpha}}= c t^{-1+\delta_1},
$$
where
\begin{equation}\label{del1}
\delta_1:= 1+ \frac{\kappa-1}{\alpha}.
\end{equation}
 Note that since $\alpha\in (1,2)$ we have $\alpha+\kappa-1>0$, which implies $\delta_1>0$\normal.
Thus,
\begin{equation}\label{lj11}
|J_1|\leq  c_3 t^{-1+\delta_1}\big(  g_{t,\theta}   * P_{t,\kappa}\big)(y-x).
\end{equation}

To estimate  $J_2$ recall that by the Taylor expansion we have
\begin{equation}
p_t^0 (x+u,y)-p_t^0(x,y)-u \cdot \nabla p_t^0(x,y)= \sum_{1\leq i,j\leq n}
u_iu_j \int_0^1 (1-\vartheta)\frac{\partial^2}{\partial
x_i\partial x_j} p_t^0(x+\vartheta u,y)d\vartheta.
\label{tay}
\end{equation}
Using \eqref{der-est}, \eqref{tay}, and the definition of $f_t$   we derive the
estimates
\begin{equation}
\begin{split}
 \Big|p_t^0(x+u,y)-p_t^0(x,y)- u\cdot \nabla  p_t^0(x,y) \Big|& \leq c_1 \Big|\sum_{1\leq i,j\leq n} u_iu_j \Big| \rho_t^2 \int_0^1 (1-\vartheta) f_t(y-x-\vartheta u)d\vartheta\\
 & \leq  n^2 \normal  c_1\|u\|^2 \rho_t^2  \int_0^1 f_t(y-x-\vartheta u)d\vartheta \\
& \leq  n^2\normal c_1 \|u\|^2 \rho_t^2 \left[ \int_0^1 e^{c_2 \vartheta\rho_t \|u\|}d\vartheta \right]  f_t(y-x)\normal
 \\&
 \leq  c_3\|u\|^2 \rho_t^2    f_t(y-x),\normal
  \end{split}
\end{equation}
where to get the last line we used that in $J_2$ we have  $\rho_t \|u\|\leq 1$.
Observe, that   for any $r>0$
\begin{align*}
\int_{r\|u\|\leq 1} (r \|u\|)^2\mu(du)&= \int_0^1 \mu\{ u: v \leq (r \|u\|)^2 \leq 1\} dv\\
&\leq  \sum_{i=1}^n \int_0^1 \mu\{ u: \, v/n \leq |ru_i| ^2\leq 1\} dv\\
& \leq n^2 \max_{1\leq i\leq n} \int_0^{1/n} \mu\{ u: \, z \leq |ru_i|^2 \leq 1 \} dz\\
&  \leq n^2 \max_{1\leq i\leq n} \int_0^1 \mu\{ u: \, z \leq |ru_i|^2 \leq 1\} dz \\
&= n^2 \max_{1\leq i\leq n}   q^L( r \ell_i) \leq n^2 \max_{1\leq i\leq n}  q^U (r \ell_i) \\
& \leq n^2 q^*(r),
\end{align*}
where  $\ell_i:=(0,\ldots, \underset{i}{1}, \ldots0) \in \mathbb{S}^n$.

Thus, using \textbf{A3},  the above calculation and that  $q^*(\rho_t)=1/t$, we can estimate  $J_2$:
\begin{align*}
|J_2|&\leq c_1 (\|y-x\|^\gamma\wedge 1) f_t(y-x)   \int_{\rho_t\| u\|\leq 1} (\|u \| \rho_t)^2 \mu(du)\\
&\leq  c_2 (\|y-x\|^\gamma \wedge 1)  f_t(y-x)q^*(\rho_t)\\
& \leq c_2\normal t^{-1}  (\|y-x\|^\kappa \wedge 1)   f_t(y-x),
\end{align*}
and by   Proposition~\ref{lem1-1}  and \eqref{growth3}   we get
\begin{equation}\label{j2}
|J_2|\leq c_3 t^{-1+\delta_2}\big(  g_{t,\theta}   * P_{t,\kappa}\big)(y-x),
\end{equation}
where  $\delta_2:=\kappa/\sigma$.

Let us estimate $J_3$. We have
\begin{equation}
\begin{split}
|J_3|&\leq  \int_{ \rho_t\|u\|>1} \big( p_t^0(x+u,y)+ p_t^0(x,y)\big) | m(y,u)-m(x,u)|\mu(du)
\\&
\quad +
 \Big|\int_{1/\rho_t <\|u\|<1}  u \cdot \nabla   p_t^0(x,y)\, [m(y,u)-m(x,u)] \mu(du)\Big|
\\&
=: J_{31}+J_{32}.
\end{split}
\end{equation}
 For $J_{31} $ we get  by \textbf{A3}, \eqref{der-est} and    Proposition~\ref{lem1-1}  the estimates
\begin{align*}
J_{31}
& \leq b_4 t^{-1}(\|y-x\|^\gamma \wedge 1) \Big\{ \int_\rn p_t^0(x+u,y)\Lambda_t(du)+ p_t^0(x,y)\Lambda_t(\rn)\Big\}\\
&\leq c_1t^{-1} (\|y-x\|^\gamma\wedge 1)  \big( (f_t *\Lambda_t)(y-x)+   n^2   f_t(y-x)\big)\\
&\leq c_2 t^{-1}  (\|y-x\|^\kappa \wedge 1)
   \Big[ (f_t*\Lambda_t)(y-x)+f_t(y-x)\big)\Big]\\
&\leq c_2 t^{-1}\Big\{  \int_\rn (\|y-x-u\|^\kappa\wedge 1)f_t(y-x-u)\Lambda_t(du) \\
& \quad + \int_\rn f_t(y-x-u)(\|u\|^\kappa \wedge 1) \Lambda_t(du) +  (\|x-y\|^\kappa \wedge 1) f_t(y-x)\Big\}\\
& \leq c_3 t^{-1+\delta_2} \Big\{\big(g_{t,\theta} *  P_{t,\kappa}*\Lambda_t\big)(y-x) + \big(f_t* P_{t,\kappa}\big)(y-x)+ \big(g_{t,\theta} * P_{t,\kappa}\big)(y-x)\Big\}\\
 &\leq   c_4  t^{-1+\delta_2}  \Big(  g_{t,\theta}   * \big(  P_t * P_{t,\kappa}+  P_{t,\kappa}\big)\Big)(y-x),\normal
\end{align*}
 where in the last line we used that $f_t(x)\leq c (g_{t,\theta}* P_t)(x)$, and that $\Lambda_t(du)$ is dominated by $P_t(du)$. \normal

Finally, we estimate $J_{32}$. By \textbf{A4}, $J_{32}=0$ in the case $\alpha\in (0,1]$.   Assume that $\alpha\in (1,2)$.  In this case, by  \textbf{A3}, \eqref{der-est}, Proposition~\ref{cor11} with $\lambda=1$,  and Proposition~\ref{lem1-1}, we get
\begin{align*}
|J_{32}|&\leq c_1 t^{-1}  (\|y-x\|^\gamma\wedge 1) \rho_t  f_t(y-x)\Big(  \int_\rn \big(\|u\|\wedge 1\big) \Lambda_t(du)\Big)\\
&\leq c_2 t^{-1}   (\|y-x\|^\kappa\wedge 1)  f_t(y-x)\\
&
\leq c_3 t^{-1}\rho_t^{-\kappa}\big(  g_{t,\theta}   * P_{t,\kappa}\big)(y-x)\\
&\leq c_3 t^{-1+\delta_2}\big(  g_{t,\theta}   * P_{t,\kappa}\big)(y-x),
\end{align*}
where in the last line we used \eqref{growth3}.

Thus, we arrive at
\begin{equation} \label{g3}
J_3 \leq c_4 t^{-1+\delta_2}   \big(  g_{t,\theta}   *  G_t\normal\big)(y-x).
\end{equation}
  Put
$$
\eta:= \delta_1\wedge \delta_2 = \frac{\kappa}{\sigma} \wedge \Big( 1+ \frac{\kappa-1}{\alpha}\Big).
$$

Thus,  combining the estimates for $J_1$, $J_2$ and $J_3$, we arrive at \eqref{ld11} with some constant    $C>0$,   $\eta$,    and $b=\theta a_2$, where $\theta\in (0,1)$ is arbitrary, and $a_2 $ is the constant from Proposition~\ref{aux1}  (cf. Remark~\ref{ak}).

\end{proof}

\subsection{Generic calculation}\label{gen-cal}

Let us rewrite the statement of Lemma~\ref{Phi-up} a bit differently. Although it might be seen as just some  technical modification, it will become clear later that this new form  allows us to write  the estimate in a rather transparent way.

Put
\begin{equation}\label{del}
\delta:= \eta/2,
\end{equation}
and
\begin{equation}\label{g20}
\tilde{g}_{t} (x) := t^{\delta} g_{t}(x).
\end{equation}
Then the estimate  \eqref{ld11} can be written as
\begin{equation}\label{Phi-up20}
\big|   \Phi_t(x,y)\big|\leq  C t^{-1+\delta}\big( \tilde{g}_t  *
G_{t}\big)(y-x),\quad   t\in (0,T], \quad x,y\in \rn.
\end{equation}

The next important step  is to estimate iteratively the convolution powers $\Phi^{\star k}$, and this is the place  where we encounter essential new difficulties. Below we explain this problem  in details, and give the generic calculation which allows us to overcome these difficulties.

Denote
\begin{equation}\label{Phi-H}
H_t(x,y):= \big(\tilde{g}_t  *
G_{t}\big)(y-x).
\end{equation}
Observe that  if this kernel would
   satisfy the following \emph{sub-convolution property}
\begin{equation}\label{H1}
\big(H_{t-s} * H_s \big)(x,y)\leq c H_t(x,y), \quad 0<s<t, \quad x,y\in \rn,
\end{equation}
then the  iterative estimation of the convolution powers $\Phi^{\star k}$ would be simple. For example, this is true  for a perturbed $\alpha$-stable noise:  in this case  we have
$$
\big|   \Phi_t(x,y)\big|\leq  C t^{-1+\delta} H_t(x,y), \quad t\in (0,T], \quad x,y\in \rn,
$$
with
 $$
H_t(x,y)=  t^{-n/\alpha}   \big( 1+ \|y-x\|/t^{1/\alpha}\big)^{ -n-\alpha},
 $$
 see  \cite{Ko89}, \cite{Ko00} and \cite{BJ07};  see also \cite{KK14a}  for more involved kernels which appear for the gradient perturbations of an $\alpha$-stable noise with $\alpha<1$.

  In our situation the kernel $H_t(x,y)$ has the more complicated structure: it is formed by the convolution of the function $g_t$ and some measure $G_t$, which seems to be inevitable because of the ``compound kernel'' structure of the  first approximation $p_t^0(x,y)$ given by (\ref{pto}).
  Moreover,  in this case we cannot in general  expect \eqref{H1} to hold true.

   To show what is going on, we  give a calculation of
   the upper bound for the convolution of two
   ``compound kernels''.

\begin{lemma}\label{phi12} Suppose that the functions  $\Phi^{i}$, $i=1,2$, satisfy for any $T>0$ the inequalities
\begin{equation}\label{Phi-up200}
\big|   \Phi_t^i(x,y)\big|\leq  C_i  t^{-1+\delta_i}\normal\big( h_t^i  *
G_{t}^i\big)(y-x),\quad t\in (0,T], \quad x,y\in \rn.
\end{equation}
with  constants $C_i>0$,\, $\delta_i>0$,\,\normal some nonnegative and integrable functions $h^i$, and some sub-probability measures  $G^i_t$,  respectively. Then $\mathfrak{F}:=\Phi^1\star\Phi^2$ satisfies
\begin{equation}\label{Phi-up21}
\big|   \mathfrak{F}_t(x,y)\big|\leq  C t^{-1+\delta}\big( \mathfrak{h}_t  *
\mathfrak{G}_{t}\big)(y-x), \quad t\in (0,T], \, x,y\in \rn,
\end{equation} with
$$
 \mathfrak{h}_t(x)\normal=\sup_{s<t}\big( h_{t-s}^1 * {h}_s^2\big) (x),
$$
$$
\delta=\delta_1+\delta_2,\quad C =C_1 C_2B(\delta_1, \delta_2),
$$
where $B(\cdot, \cdot)$ is the Beta-function,
and
$$
\mathfrak{G}_t(dw):= {1\over B(\delta_1, \delta_2)}\int_0^1 \int_\rn (1-r)^{-1+\delta_1} r^{-1+\delta_2}  G_{t(1-r)}^{1}(dw-u) G_{tr}^2(du)dr.
$$
Moreover,   $\mathfrak{G}_t(dw)$ is the  sub-probability  measure, i.e. $\mathfrak{G}_t(\rn) \leq 1$ for each $t\in [0,T]$.

\end{lemma}
\begin{proof}  Making the change of variables, we derive
\begin{align*}
\big|\mathfrak{F}_t(x,y)\big|&\leq C_1 C_2 \int_0^t \int_{\real^{3n}} h_{t-s}^1 (z-w_1-x) h_s^2(y-w_2-z) \frac{ G_{t-s}^1 (dw_1) G_s^2 (dw_2)}{(t-s)^{1-\delta_1} s^{1-\delta_2}}dzds\\
&\leq C_1C_2 B(\delta_1,\delta_2) t^{-1+\delta_1+\delta_2} \int_\rn\mathfrak{h}_t(y-x-w) \Big[ \int_0^1 \int_\rn \frac{G_{t(1-r)}^{1}(dw-u) G_{tr}^2(du)}{B(\delta_1, \delta_2)(1-r)^{1-\delta_1} r^{1-\delta_2}} dr\Big],
\end{align*}
which gives \eqref{Phi-up21}.  Further, since   $G_t^i (\rn)\leq 1$, $i=1,2$, we  get $\mathfrak{G}_t(\rn)\leq 1$.

\end{proof}

By Lemma \ref{phi12},     we have the following estimates for the convolution powers of $\Phi$:
\begin{equation}\label{Phi10}
\big|   \Phi_t^{\star k} (x,y)\big|\leq {C^k \Gamma^k(\delta)\over   \Gamma(k\delta)  } t^{-1+\delta k} \big(  h_{t}^{(k)}   *
{G}_{t}^{(k)}\big)(y-x), \quad t\in (0,T], \quad x,y\in \rn,\quad k\geq 1,
\end{equation}
where the constant $C>0$ comes from \eqref{Phi-up20},    $h^{(1)}\equiv \tilde{g}$,
$$
  h^{(k+1)}_t(x)=\sup_{s<t}\big(h^{(k)}_{t-s}*h^{(1)}_s\big)(x),\quad k\geq 1,
$$
and
\begin{equation}\label{gk-new}
    G_{t}^{(k)}(du):=
    \begin{cases}
    G_{t}(du)\normal &  k=1,\\
    {1\over B(\delta,   (k-1)  \delta)}\int_0^1 \int_\rn (1-r)^{-1+(k-1)\delta} r^{-1+\delta}  G_{t(1-r)}^{(k-1)}(dw-u) G_{tr}(du)dr,& k\geq 2,
\end{cases}
\end{equation}
 Hence, to guarantee the convergence of the series of convolution powers (\ref{Psi}), it is enough to derive a proper upper bound on the sequence of functions   $h^{(k)}_t$. At this concern, we give the following lemma.

\begin{lemma}\label{con-g-new}
Let   $g_t(x)$ be  defined in  \eqref{gt10}.   Then for any $\theta\in (0,1)$, $T>0$,   we have
  \begin{equation}\label{gt11-new}
  (g_{t-s} *g_{s} )(x) \leq C_0(\theta) g_{t}(\theta x), \quad 0<s<t, \quad x\in \rn,\, t\in (0,T],
\end{equation}
where $C_0(\theta)=c(1-\theta)^{-n}$, and $c>0$  is some constant.
\end{lemma}
\begin{proof}
Consider the integral
    \begin{equation}
    I(t,x):= \int_\rn g_{t-s}(x-z)g_{s}(z)dz. \label{g1-1}
    \end{equation}
 Suppose that $0<s\leq t/2$.  Note that for $s<\frac{t}{2}$ we have by monotonicity of $\rho_t$ that $\rho_{t-s}\leq \rho_{t/2}$. Further, for $c_1\geq 1$ we have  $q^*(r)\leq q^*(c_1r) \leq c_1^2  q^*(r)$ for all $r\geq 1$, which implies $\rho_t\asymp \rho_{c_1t}$ for all $t\in (0,T]$. Therefore,
 $$
 \rho_{t/2}\leq c_2 \rho_t, \quad t\in (0,T].
 $$
Since  $\rho_t$ is   decreasing,      the triangle inequality
  $$
    \|x-z\|\rho_{t-s} +\|z\|\rho_s \geq \|x\| \rho_t
    $$
      gives
    \begin{equation}
    \begin{split}
    I(t,x)&\leq  e^{-b\theta\|x\|\rho_t}
   \int_\rn   \rho_{t-s}^n \rho_s^n e^{-b (1-\theta)[\rho_{t-s}\|x-z\|+\rho_s \|z\|]}dz
     \\
     &\leq    \rho_{t/2}^n  e^{-b \theta\rho_t\|x\|}   \int_\rn \rho_s^n e^{-b (1-\theta)\rho_s \|z\|}dz\\
     &= C(\theta) g_{t}(\theta x),\label{g22}
    \end{split}
    \end{equation}
  where $C(\theta):=  c_2^n  c_0  [b(1-\theta)]^{-n}$, $c_0:= \int_\rn e^{-\|z\|}dz$.
\end{proof}

\begin{lemma}\label{Phi13} For any $T>0$ and  any sequence $(\theta_k)_{k\geq 1}$  such that
$$
\theta_1=1, \quad \text{and} \quad \theta_{k+1}<\theta_{k}, \quad \theta_k>0, \quad k\geq 1,
$$
 one has for $k\geq 2$ the estimate
\begin{equation}\label{Phi10a}
\big|   \Phi_t^{\star k} (x,y)\big|\leq C_k t^{-1+\delta k} \big(g_{t}^{(k)}*
{G}_{t}^{(k)}\big)(y-x), \quad t\in (0,T], \quad x,y\in \rn,
\end{equation}
where the   sub-probability    measures $(G^{(k)})_{k\geq 1}$ are defined in \eqref{gk-new},
\begin{equation}\label{gt-new}
  g_t^{(k)}(x):=t^{\delta k } g_t(\theta_k x),\quad k\geq 1,
\end{equation}
and
$$
 C_k:={C^k c^{k-1} \Gamma^k(\delta)\over \Gamma(k \delta)}\prod_{j=2}^k \left(\frac{1}{\theta_{j-1}-\theta_j}\right)^n, \quad k\geq 2,
$$
where the positive constants $C$  and $c$ come from  \eqref{Phi-up20} and Lemma~\ref{con-g-new}, respectively.
\end{lemma}
\begin{proof}
  By monotonicity of $g_t(x)$ in $x$ we have $g_t(x)\leq g_{t}(\theta_{k-1}x)$.
Therefore, using Lemma~\ref{con-g-new}   with $\theta=\theta_k/\theta_{k-1}$\normal\,, we get
$$
\begin{aligned}
(g_{t-s}^{(k-1)}*g_s^{(1)})(x)
&  \leq \normal t^{\delta k}  \int_{\rn}g_{t-s}(\theta_{k-1}x-\theta_{k-1}y) g_{s}(\theta_{k-1}y)\, dy\\
&=   t^{\delta k}   \theta_{k-1}^{-n}\int_{\rn}g_{t-s}(\theta_{k-1}x-y') g_{s}(y')\, dy'\\
&  \leq   D_k    t^{\delta k}    g_{t}(\theta_{k} x)\normal\\
&= D_k g_t^{(k)}(x),
\end{aligned}
$$
where
$$
D_k= (\theta_{k-1})^{-n}C_0(\theta_k/\theta_{k-1})
= \frac{c}{\big(\theta_{k-1}(1-\theta_k/\theta_{k-1})\big)^n}
=\frac{c}{\big(\theta_{k-1}-\theta_k\big)^n}.
$$
Then \eqref{Phi10a} follows from (\ref{Phi10}).
\end{proof}
 Estimate \eqref{Phi10a} is still hardly applicable for verifying  the convergence of the series of convolution powers (\ref{Psi}): to keep this sequence of estimates consistent, one should choose the sequence $\{\theta_k\}$ such that $\inf_k\theta_k>0$, and then it is difficult to bound properly the values of the constants $C_k$. In order to illustrate this, take e.g. $\theta_k:= \frac{1}{2} +\frac{1}{2k}$. Then
  $$
   \prod_{j=2}^{k} \Big(\frac{1}{\theta_{j-1}-\theta_j}\Big)^n   = \prod_{j=2}^{k}( 2j(j-1))^n=\Big( 2^k k! (k-1)! \Big)^n,
  $$
  which increases faster than $\Gamma(k\delta)$ in the denominator   in the definition of $C_k$.

 In order to overcome this problem,    we  change after finite number of steps the sequence $g_t^{(k)}(x)$. This change  finally allows to prove the convergence of the series (\ref{Psi}).

 Let
 \begin{equation}\label{k0}
 k_0=\left[ \frac{n}{\alpha \delta}\right] +1.
 \end{equation}
Note that for such $k_0$ we have  $  t^{\delta k_0}   \rho_t^n \leq   c(k_0)\normal$   for $t\in (0,T]$.
Then
  \begin{align*}
\Big(g_{t-s}^{(k_0)} * g_s^{(1)}\Big)(x)&   \leq  c(k_0) \int_\rn s^\delta \rho_s^n  e^{-b \theta_{k_0} \rho_{t-s} \|x-z\| -b \rho_s \|z\| } dz\\
 & \leq  c(k_0)  M e^{-b\zeta \rho_t\|x\|}   =  c(k_0)M  \rho_t^{-n} g_{t,\zeta}(x),\quad 0<s<t,
 \end{align*}
 where  $\zeta=\theta_{k_0}$, $g_{t,\zeta}(x)$ is of the form  \eqref{gt},\,\normal
 and
 \begin{equation}\label{M}
 M:=   T^\delta \normal \int_\rn e^{-b(1-\zeta) \|z\|}dz.
 \end{equation}
 By induction, we get
 \begin{equation}\label{gtk0-new}
  \bar{g}_t^{(k_0+\ell+1)}  (x)  := \sup_{0<s<t}\Big(\bar{g}_{t-s}^{(k_0+\ell)} * g_s^{(1)}\Big)(x)\normal
\leq c(k_0)  M^{\ell+1}    \rho_t^{-n}  g_{t,\zeta}(x),  \quad  \ell\geq 0.
\end{equation}

\begin{lemma}\label{gtkl}
For any $T>0$ we have
\begin{equation}\label{LZkl}
\big| \Phi_t^{\star (k_0+\ell)}(x,y) \big|\leq  D_\ell  t^{-1+ \delta (k_0+\ell)}  \rho_t^{-n} \big( g_{t,\zeta} * G_t^{(k_0+\ell)} \big) (y-x), \quad \ell\geq 1,\,  t\in (0,T],
\, x,y\in \rn,\normal
\end{equation}
where $k_0$ is given by  \eqref{k0}, the family of   sub-probability    measures $\{ G_t^{(k)},\, t>0, \, k\geq 1\}$ is defined in \eqref{gk-new},
\begin{equation}\label{Dl}
 D_\ell := \frac{ C(k_0) (CM)^\ell  \Gamma^{k_0+\ell}(\delta)}{\Gamma((k_0+\ell)\delta)},\normal
\end{equation}
where $C(k_0)>0$ is some constant, and  $C,M>0$ come, respectively, from \eqref{Phi-up20} and \eqref{M}.
\end{lemma}
The proof follows by induction; we omit the details.

\subsection{Proof of Theorem~\ref{t-main1}}\label{pr-main}

 a) Using \eqref{Phi10a}, \eqref{LZkl} and that   $g_t^{(k)}(x)\leq  T^{k\delta}  g_{t,\zeta}(x)$\normal, $1\leq k\leq k_0$, we get   for $t\in (0,T]$\normal
\begin{equation} \label{Psi-est}
\begin{split}
|\Psi_t(x,y)|&\leq \sum_{k=1}^\infty \big|\Phi_t^{\star k} (x,y)\big|
\leq \sum_{k=1}^{k_0} C_k t^{-1+k\delta} \big( g_t^{(k)}* G_t^{(k)}\big) (y-x)\\
&\quad + \sum_{\ell=1}^\infty  D_\ell t^{-1+\delta (k_0+\ell) } \rho_t^{-n} \big( g_{t,\zeta} * G_t^{(k_0+\ell)}\big) (y-x)\\
&\leq t^{-1+\delta}   \Big( g_{t,\zeta} * \Big(\sum_{k=1}^{k_0} T^{\delta(k-1)} C_k G_t^{(k)}+ \sum_{\ell=1}^\infty  T^{\delta(k_0+\ell-1)}  D_\ell G_t^{(k_0+\ell)}\Big) \Big) (y-x)\normal \\
&\leq   \big( \sum_{k=1}^{k_0} T^{\delta(k-1)} C_k + \sum_{\ell=1}^\infty T^{\delta(k_0+\ell-1)} D_\ell\big)\normal  t^{-1+\delta} \big(g_{t,\zeta}* \Pi_t\big)(y-x),
\end{split}
\end{equation}
\begin{equation}\label{Pi}
\Pi_t(du):=  \frac{   \sum_{k=1}^{k_0} T^{\delta(k-1)}C_k G_t^{(k)}(du)+ \sum_{\ell=1}^\infty T^{\delta(k_0+\ell-1)}  D_\ell G_t^{(k_0+\ell)} (du)}{ \sum_{k=1}^{k_0} T^{\delta(k-1)} C_k + \sum_{\ell=1}^\infty T^{\delta(k_0+\ell-1)} D_\ell }, \quad t\in (0,T].\normal
\end{equation}
Since $G_t^{(k)}$, $k\geq1$, are the sub-probability measures, then
\begin{equation}
\Pi_t(\rn)  \leq   1, \quad t\in (0,T]. \label{pi2}
\end{equation}
 Thus, we proved that the series $\Psi_t(x,y)= \sum_{k=1}^\infty \Phi_t^{\star k} (x,y)$ converges for any   $t\in (0,T]$ ,   $x,y\in \rn$, uniformly on compact subsets  of $(0,\infty) \times \rn\times \rn$.

Finally, let us show that  $(p^0 \star \Psi)_t(x,y)$ is well defined.
Using the upper bound for $p^0$ (cf. \eqref{der-est} with $k=0$) and  the estimate for $\Psi$ (cf. \eqref{Psi-est}), we get by Lemmas~\ref{phi12} and \ref{con-g-new}   \normal
\begin{equation}
\big|(p^0 \star \Psi)_t(x,y)\big|\leq  t^\delta \big(g_{t,\chi} * \tilde{\Pi}_t\big) (y-x), \quad t\in (0,T], \, x,y\in \rn, \label{up-11}
\end{equation}
where   $\chi\in (0,\zeta)$\normal \, ($\zeta$ comes from \eqref{Psi-est}, and    $g_{t,\kappa}(x)$ is of the form \eqref{gt},\normal\, and
\begin{equation}\label{til-pi}
\tilde{\Pi}_t(dw):= \delta \normal \int_0^1 \int_\rn (1-r)^{-1+\delta} \Pi_{t(1-r)}(dw-u)P_{tr}(du)dr.
\end{equation}
 By the definition  of $\Pi_t$ and  $P_t$, $\tilde{\Pi}_t(\rn)\leq 1 $.   Thus, the expression \eqref{sol10} is well defined, and the series involved in this expression converge  absolutely, uniformly on compact sets of $(0,\infty)\times \rn\times\rn$.  \normal

b)    By assumption \textbf{A2} we have
$$
\|\xi\|^k\big|  e^{i\xi (x-y)-t q(y,\xi)} \big| \leq e^{- c t_0 q(\xi)}, \quad x,y\in \rn,\,  t\in [t_0,\infty),
$$
for any $t_0>0$  and $k\geq 0$. Therefore, since the function $|e^{i\xi (x-y)-t q(y,\xi)}|$ is   continuous,  has continuous derivatives,   and $q(\xi)\geq c \|\xi\|^\alpha$ for  $\|\xi\|\geq 1$ (see \eqref{psipm1} and \eqref{qal}), then by the dominated convergence theorem the function
$$
p_t^0(x,y)= (2\pi)^{-n} \int_\rn e^{i\xi (x-y) - tq(y,\xi)}d\xi \normal
$$
has   continuous derivatives in $(t,x)\in (0, \infty)\times \rn$.

Next we prove that the functions $\Phi_t^{\star k}(x,y)$, $k\geq 0$,  are continuous in $(t,x,y)$ on $(0, \infty)\times \rn\times \rn $.
 Let us show that $\Phi_t(x,y)$ is continuous. As we just have shown, for any $t_0>0$
 $$
 \Big| \partial_{x_i x_j}^2 p^0_t(x,y)\Big| \leq c, \quad t\geq t_0, \, x,y\in \rn, \, 1\leq i,j\leq n,
 $$
 then
 $$
   \Big|  p_t^0(x+u,y)-  p_t^0(x,y)-\nabla_x  p^0_t(x,y) \cdot u \I_{\{\|u\|\leq 1\}} \Big|\normal \leq C\big(\|u\|^2\wedge 1 \big), \quad t>t_0, \, x,y\in \rn
 $$
Therefore, by the dominated convergence theorem we derive that $L_x p^0_t(x,y)$ is continuous in $(t,x,y)$ on $[t_0, \infty)\times \rn\times \rn $. Therefore, since $\partial_tp_t^0(x,y)$ is continuous in $(t,x,y)$ on $[t_0, \infty)\times \rn\times \rn $ and $t_0>0$ is arbitrary,      we  derive \normal \, the desired  continuity of $\Phi_t(x,y)$.

To show that the convolutions $\Phi_t^{\star k}(x,y)$ are continuous, we use induction.

Suppose that $\Phi_t^{\star (k-1)}(x,y)$ is continuous.
Let $t_0>\eps>0$,  and suppose that $t\in [t_0,\infty)$.  Write
\begin{equation}\label{Phik10}
\begin{split}
\Phi_t^{\star k}(x,y)&=\int_0^{t-\eps} \left[\int_\rn \Phi_{t-s}^{\star (k-1)} (x,z)\Phi_s(z,y)dz\right]ds\\
&\quad +\int_{t-\eps}^t \left[\int_\rn \Phi_{t-s}^{\star (k-1)} (x,z)\Phi_s(z,y)dz\right]ds\\
&=\int_0^{t-\eps} \left[\int_\rn \Phi_{t-s}^{\star (k-1)} (x,y-z)\Phi_s(y-z,y)dz\right]ds\\
&\quad +\int_0^{\eps} \left[\int_\rn \Phi_{s}^{\star (k-1)} (x,x-z)\Phi_{t-s}(x-z,y)dz\right]ds\\
&=I_1(t,x,y)+I_2(t,x,y).
\end{split}
\end{equation}
We prove the continuity of $I_1(t,x,y)$, the continuity of $I_2(t,x,y)$ follows by the same argument.

By the induction assumption, the function $\Phi_{t-s}^{\star (k-1)} (x,z)\Phi_s(z,y)$ is continuous in $s\in (0,t-\eps]$, $t\in [t_0,\infty)$, $(x,y)\in  \rn\times \rn$.   Moreover,
by \eqref{Phi10a}   and \ref{ld11}\normal
\begin{align*}
|\Phi_{t-s}^{\star (k-1)}(x,y-w) \Phi_s(y-w,y)|&\leq C   (t-s)^{-1+2\delta (k-1)} s^{-1+\eta}\rho_{t-s}^n  \big( g_s *G_s\big)(w)\normal \\
&\leq C(\eps)   s^{-1+\eta} \big( g_s *G_s\big)(w)\normal.
\end{align*}
Since the right-hand side of the above inequality is integrable on $[0,t-\eps]\times \rn$, we get by the dominated convergence theorem that $I_1(t,x,y)$ is continuous in $(t,x,y)\in [t_0,\infty)\times \rn\times \rn$. Finally, since $t_0$ and $\eps$ are  arbitrary, we get the continuity in $(t,x,y)$ on $(0,\infty)\times \rn\times\rn$.

Since the series $\sum_{k=1}^\infty \Phi_t^{\star k}(x,y)$ converges uniformly on compact subsets of $(0,\infty)\times \rn\times \rn$, the function $\Psi_t(x,y)$ is continuous in $(t,x,y)\subset (0,\infty)\times \rn\times \rn$.

The proof of the continuity of $p^0 \star \Psi $ follows by the same argument as the proof of continuity of $\Phi_t^{\star k}$, we only need to use estimates on $p^0$ and $\Psi$, see  \eqref{der-est} and \eqref{Psi-est}).  \normal

 \qed

\subsection{Continuity properties of the operator $S_t$}\label{pr-cont}

Note that by construction  we have  for any sub-probability measure $M_t(\cdot)$ on $\rn$
\begin{equation}\label{int-est}
\int_\rn \big(g_{\theta, t}* M_t\big) (y)dy \leq C, \quad t\in (0,T],\,  x\in \rn.
\end{equation}
Therefore, by \eqref{der-est}  and \eqref{up-11}  we have
$$
\int_\rn p_t(x,y)dy\leq C, \quad t\in (0,T], \, x\in \rn.
$$
Then  the operator $S_tf $, $t>0$,  (cf.  \eqref{t-semi}) is well defined for  any bounded measurable function $f$.

\begin{lemma}\label{P-cont}
\begin{enumerate}\label{cont}
  \item For any $t>0$ the operator  $S_t$ maps $ C_\infty(\rn)$ into $C_\infty(\rn)$.

  \item For every $f\in C_\infty(\rn)$  we have $\lim_{t\to 0+} \|S_t f-f\|_\infty=0$.

\end{enumerate}
\end{lemma}
The proof  relies on the proposition below.
\begin{proposition}\label{H-bound-a}
For every $f\in C_\infty(\rd)$
\begin{equation}\label{Fel}
\lim_{|x|\to\infty} \int_\rd p_t^0(x,y)f(y)dy =0, \quad \text{for any}\quad  t>0,
\end{equation}
\begin{equation}\label{Fel2}
\sup_{x\in \rd} \Big|\int_\rd p_t^0(x,y)f(y)dy -f(x)\Big|\to 0, \quad t\to 0.
\end{equation}
\end{proposition}
In order to keep the presentation as clear as possible, we defer the proof of this Proposition to Appendix B.
\begin{proof}[Proof of Lemma~\ref{P-cont}]
1. The continuity of $S_tf$ follows from the  continuity of $p_t(x,y)$.   To prove that $S_tf(x)$ vanishes as $\|x\|\to \infty$, we use  the representation for $p_t(x,y)$ (cf. \eqref{sol10} and \eqref{r}):
 \begin{align*}
\Big| \int_\rn p_t(x,y)f(y)dy \Big| & \leq \int_\rn  p_t^0(x,y) |f(y)|dy + \int_\rn  |(p^0\star \Psi)_t(x,y) f(y)|dy\\
&= I_1(t,x)+I_2(t,x).
\end{align*}
By the first statement of Proposition~\ref{H-bound-a}, the first term on the right-hand side tends to 0 as $\|x\|\to\infty$.
Using  the upper estimate on  $p^0\star \Psi$ (cf. \eqref{up-11}) we get
\begin{align*}
I_2(t,x)&\leq C t^\delta \int_\rn \Big[ \int_\rn g_{t,\chi}(y-x-w)|f(y)|dy\Big] \tilde{\Pi}_t(dw)\\
&= C t^\delta \int_\rn \Big[ \int_\rn g_{t,\chi}(z-w)|f(z+x)|dy\Big] \tilde{\Pi}_t(dw), \end{align*}
and the right-hand side tends to zero as $|x|\to\infty$ by the dominated convergence theorem (recall that the parameter $\chi$ comes from \eqref{up-11}).

2.  By \eqref{Fel2} it is enough to show that
\begin{equation}\label{Fel3}
\sup_x \Big|\int_0^t \int_\rd \int_\rd \Psi_{t-s}(x,z)p_s^0(z,y)f(y)dydzds\Big| \to0, \quad t\to 0.
\end{equation}
By \eqref{up-11} we have
\begin{align*}
\sup_x \Big|\int_0^t \int_\rd \int_\rd \Psi_{t-s}(x,z)p_s^0(z,y)f(y)dydzds\Big|&\leq
 c_1t^\delta \int_\rn \big(g_{t,\chi} * \tilde{\Pi}_t\big) (y)\,dy\leq c_2 t^\delta,
\end{align*}
which finishes the proof.

\end{proof}

 \normal

\section{Approximative positive maximum principle. Proof of Theorem~\ref{t2}. }\label{ver}

We follow with the necessary changes the approach described in  \cite{KK14a}.

So far in Section~\ref{constr} we constructed the function $p_t(x,y)$ under the assumption that $p_t(x,y)$ is a fundamental solution to the Cauchy problem for $\prt_t -L$. As in \cite{KK14a}, the straightforward way  to check \eqref{L_delta} and \eqref{L_fund} meets difficulties. To explain them, let us look at the begaviour of the derivatives of $p_t^0(x,y)$ near the origin.

Let
\begin{equation}\label{mathp}
\mathcal{P}_t (du):= c\Big(P_t(du)+  (P_t *\Lambda_t)(du)\Big).
\end{equation}
where  for any $T>0$   the constant $c>0$ is chosen such that   $\mathcal{P}_t(\rn)\leq 1$\,   for all $t\in(0,T]$,  which is possible since $\Lambda_t(\rn)\leq n^2$ (cf.  Proposition~\ref{lem1}). \normal

\begin{proposition}\label{p0-der}
The function $p_t^0(x,y)$ is differentiable with respect to $t$, the derivative $\partial_t p_t^0(x,y)$ is  continuous in  $(t,x,y)\in(0,\infty)\times \rn\times \rn$, and for $k=k_1+\dots+k_n\geq 0$,  $T>0$, \normal  there exist constants $\tilde{A}_k, \tilde{a}_k>0$ such that
\begin{equation}\label{p0der-es}
\Big|\frac{\partial }{\partial t} \frac{\partial^k}{\partial x_1^{k_1} \ldots \partial x_n^{k_n}}   p_t^0(x,y)\Big|\leq t^{-1}\rho_t^k \big(f_{up} * \mathcal{P}_t\big)(y-x), \quad   t\in (0,T]\normal, \quad x,y\in \rn,
\end{equation}
where $\mathcal{P}_t(du)$ is the   sub-probability measure  \,  defined in \eqref{mathp}, and $f_{up}$ is of the form \eqref{flu} with  constants $\tilde{A}_k$ and $\tilde{a}_k$ in the place of $d_1$ and $d_2$, respectively.
\end{proposition}
The proof of this proposition can be obtained by a simple modification of the proof of the respective statement in \cite{KK12b}, see also \cite{K13}. We omit the details.

From  \eqref{p0der-es} and \eqref{Psi-est}   it is unclear why the function $p_t(x,y)$ given by \eqref{sol10} is in the domain of the operator $L$, which was so far defined on $C_\infty^2(\rn)$-functions. Consequently, we cannot check straightforwardly that \eqref{L_fund} holds true,  and in such a way to verify the correctness of the procedure performed in Section~\ref{constr}.

  To avoid this difficulty, we introduce for $\eps>0$  the  auxiliary function
\begin{equation}\label{pe}
p_{t,\epsilon}(x,y):=p_{t+\epsilon}^0(x,y)  + \int_0^t \int_{\rn}  p_{t-s+\eps}^0(x,z)  \Psi_s(z,y) dzds.
\end{equation}
Since the  additional time shift by positive $\eps$ removes the singularity at the point $s=t$, we have the following properties of $p_{t,\eps}(x,y)$:

\begin{itemize}

\item[(i)]   $p_{\cdot,\epsilon}(x,y)\in C^1((0,\infty))$ for any fixed $\epsilon>0$, $x,y\in \rn$;

\item[(ii)]  $p_{t,\epsilon}(\cdot,y)\in C^2_\infty(\rn)$ for any fixed $\epsilon>0$, $t>0$, $y\in \rn$;

\item[(iii)] for any $0<\tau<T$ we have $ p_{t,\epsilon}(x,y)\to p_t(x,y)$  as  $ \epsilon\to 0$, uniformly in
 $(t,x)\in [\tau,T]\times \rn\times \rn$.

\item[(iv)] for any $0<\tau<T$ we have
 $$
 q_{t,\epsilon}(x,y):=\big(\prt_t-L_x\big) p_{t,\epsilon}(x,y)\to 0,  \quad \eps\to 0,
 $$
  uniformly in  $(t,x,y)\in [\tau,T]\times \rn\times \rn $.
\end{itemize}

The proofs of the above properties are completely analogous to the proofs of the respective properties for
\begin{equation} \label{Pte}
S_{t,\eps} f(x):=\int_\rn p_{t,\epsilon}(x,y)f(y)dy, \quad t>0, \, x\in \rn,
\end{equation}
given in Lemmas~\ref{aux-ep} and \ref{l5}. Here we only mention that properties iii) and iv) motivate the name \emph{ the approximate fundamental solution} used for $p_{t,\eps}(x,y)$.

\begin{lemma}\label{aux-ep}
\begin{enumerate}
\item  For any  $f\in C_\infty(\rn)$, $\eps>0$, the function $S_{t,\eps}f(x)$ belongs to $C^1((0,\infty))$ as a function of $t$, and
to $C^2_\infty(\rn)$ as a function of $x$.
 \item    For every   $f\in C_\infty(\rn)$, $T>0$,
  \begin{equation}\label{conv_pte}
 \lim_{\eps\to 0} \|S_{t,\eps}f- S_t f\|_\infty=0,
  \end{equation}
   uniformly in $t\in [0,T ]$, and for every $\eps>0$
 \begin{equation}\label{conv_x}
  S_{t,\eps}f(x)\to 0, \quad \|x\|\to \infty,
\end{equation}
 uniformly in $t\in [0,T]$.

 \item For $f\in C_\infty(\rn)$ we have
 $$
 \lim_{t,\eps\to 0+} \|S_{t,\eps} f-f\|_\infty =0.
  $$
 \end{enumerate}
\end{lemma}
\begin{proof}
The proof of the first statement follows from the upper estimate \eqref{Psi-est} on $\Psi_t(x,y)$, Proposition~\ref{p0-der}, and the dominated convergence theorem.

Observe that the function
$$
[0,T]\ni t\mapsto \int_{\rn}p_{t}^0(\cdot,y)f(y)\, dy\in C_\infty(\rn)
$$
is continuous, since the function $p_{t}^0(x,y)$  is continuous in $t$ for $t>0$, and the continuity of the integral at $t=0$ is provided by Proposition~\ref{H-bound-a}.  Then
$$
 \int_{\rn}p_{t+\epsilon}^0(x,y)f(y)\, dy\to \int_{\rn}p_{t}^0(x,y)f(y)\, dy, \quad \epsilon\to 0,
$$
uniformly in $t\in [0,T], x\in \rn$. This together with  Proposition~\ref{p0-der}, estimate \eqref{p0der-es}  and the dominated convergence theorem implies statement 2.

 The proof of statement 3
 is  a slight variation of the proof of statement 2 in Lemma \ref{cont}: we just need to substitute $p_{t}^0(x,y)$ and $\Psi_{t-s}(x,z)$ with $p_{t+\eps}^0(x,y)$ and $\Psi_{t+\eps-s}(x,z)$  in \eqref{Fel} and \eqref{Fel3}, respectively.

\end{proof}

Denote
\begin{equation}\label{qte20}
V_{t,\eps} f(x)= \big(\partial_t-L_x\big) S_{t,\eps} f (x), \quad f\in C_\infty(\rn).
\end{equation}

\begin{lemma}\label{l5}  For any  $f\in C_\infty(\rn)$  we have
\begin{enumerate}

  \item
  \begin{equation}\label{conv_loc}
V_{t,\eps} f(x)\to 0, \quad \epsilon\to 0,
\end{equation}
uniformly in  $(t,x)\in [\tau,T]\times \rn$ for any $\tau>0$, $T>\tau$;

\item
 \begin{equation}\label{conv_int}
\int_0^t V_{s,\eps} f(x)ds\to 0, \quad \epsilon\to 0,
\end{equation}
uniformly in  $(t,x)\in [0,T]\times \rn$ for any $T>0$.
\end{enumerate}
\end{lemma}
\begin{proof}  Note that $S_{t,\eps} f\in C^2_\infty(\rn)$, and  thus the expression
\begin{equation}\label{LP}
LS_{t,\eps}f(x)=L_x \int_\rn  p_{t+\epsilon}^0(x,y)f(y)dy + L_x\int_0^t \int_\rn\int_\rn p_{t-s+\epsilon}^0(x,z)\Psi_s(z,y)f(y)dydzds
\end{equation}
is well defined. Let us show that we can interchange  $L_x$ with the integrals in \eqref{LP}, i.e.  that
 \begin{equation}\label{LP1}
LS_{t,\eps}f(x)=\int_\rn  L_x p_{t+\epsilon}^0(x,y)f(y)dy + \int_0^t \int_\rn\int_\rn L_x p_{t-s+\epsilon}^0(x,z)\Psi_s(z,y)f(y)dydzds.
\end{equation}
Recall the representation of $L$, cf. \eqref{gener}. By Proposition~\ref{p0-der}, the gradient term in \eqref{gener}  can be interchanged with the integral by the dominated convergence theorem.  To do the same with the ``integral part'' $\mathcal{L}$  of $L$, observe that
$$
\mathcal{ L}f(x)=\lim_{\upsilon\to 0+}\mathcal{L}^{(\upsilon)} f(x), \quad \mathcal{L}^{(\upsilon)}f(x):=\int_{\|u\|>\upsilon}(f(x+u)-f(x)-\nabla f(x) \cdot u \I_{\{\|u\|\leq 1\}}\big)\mu(x,du).
  $$
Clearly, the operator $L^{(\upsilon)}$ can be interchanged with the integrals by  the Fubini theorem. On the other hand,
  \begin{align*}
 | \mathcal{ L}f(x)-\mathcal{ L}^{( \upsilon)}f(x)|&=\left|\int_{\|u\|\leq \upsilon}\Big(f(x+u)-f(x)-\nabla f(x) \cdot u \I_{\{\|u\|\leq 1\}}\Big)\mu(x,du)\right|
 \\&\leq C\sup_{x\in \rn}\|\nabla^2 f(x)\|\int_{\|u\|\leq \upsilon}\|u\|^2 \mu(x,du)
\end{align*}
Using again Proposition~\ref{p0-der}, \eqref{Psi-est},  and the dominated convergence theorem, we can pass to the limit in
$$
\int_0^t \int_\rn\int_\rn \mathcal{ L}_x^{(\upsilon)} p_{t-s+\epsilon}^0(x,z)\Psi_s(z,y)f(y)dydzds
$$
as $\upsilon\to 0$. Thus,  \eqref{LP1} holds true.

Similarly,  by Proposition~\ref{p0-der} we get
\begin{equation}\label{TP1}
\begin{split}
\prt_t S_{t, \eps}f(x)&=\int_\rn \prt_tp_{t+\epsilon}^0(x,y)f(y)dy + \int_0^t \int_\rn\int_\rn  \prt_t p_{t-s+\epsilon}^0(x,z)\Psi_s(z,y)f(y)dydzds\\
&\quad + \int_\rn\int_\rn p_{\epsilon}^0(x,z)\Psi_t(z,y)f(y)dydz.
\end{split}
\end{equation}
Since
$$
(L_x-\prt_t)p_t^0(x,y)=\Phi_t(x,y),
$$
combining \eqref{LP1} and \eqref{TP1} we derive
\begin{equation}\label{qte10}
\begin{split}
V_{t,\epsilon}f(x)& =\int_\rn \int_\rn p_{\epsilon}^0(x,z)\Psi_t (z,y)f(y)dydz- \int_\rn \Phi_{t+\epsilon} (x,y)f(y)dy\\
&-\int_0^t \int_\rn\int_\rn  \Phi_{t-s+\epsilon}(x,z) \Psi_s(z,y)f(y)dydzds.
\end{split}
\end{equation}
Since the  function $\Psi$ satisfies the equation
$$
\Phi_t(x,y)=\Psi_t(x,y)-\int_0^t\int_{\rn}\Phi_{t-s}(x,z)\Psi_s(z,y)\,dzds,
$$
we can rewrite $V_{t,\epsilon}f(x)$ as follows:
\begin{align*}
V_{t,\epsilon}f(x)& = \int_\rn \left(\int_\rn p_{\epsilon}^0(x,z)\Psi_t (z,y)dz-\Psi_{t+\epsilon} (x,y)\right)f(y)dy\\
&\qquad +\int_\rn \left(\int_t^{t+\epsilon} \int_\rn \Phi_{t-s+\epsilon} (x,z)\Psi_s(z,y)dzds\right) f(y)dy\\&=:V_{t,\epsilon}^1f(x) +V_{t,\epsilon}^2f(x).
\end{align*}
By the uniform  continuity of  $\Psi$ on compact subsets of $(0,\infty)\times \rn\times \rn$ and estimate \eqref{Psi-est},  we have
 for $f\in C_\infty(\rn)$
$$
\sup_{t\in [\tau, T], x\in \rn}\left|\int_{\rn}\Psi_{t+\epsilon} (x,y)f(y)\, dy-\int_{\rn}\Psi_{t} (x,y)f(y)\, dy\right|\to 0, \quad \eps\to 0.
$$
Similarly,  \eqref{Psi-est}, \eqref{Fel2}, and the uniform continuity of $\Psi$ on compact subsets of  $(0,\infty)\times \rn\times \rn$ gives
$$
\sup_{t\in [\tau, T], x\in \rn\times \rn}\left|\int_{\rn} p_\eps^0(x,z)\Psi_{t} (z,y)\, dz-\Psi_{t} (x,y)\right|\to 0, \quad \eps\to 0,
$$
$$
\sup_{t\in [\tau, T], x\in \rn}\left|\int_{\rn}\int_{\rn}p_\eps^0(x,z)\Psi_{t} (z,y)f(y)\, dzdy-\int_{\rn}\Psi_{t} (x,y)f(y)\, dy\right|\to 0, \quad \eps\to 0.
$$
This proves  (\ref{conv_loc}) with   $V_{t,\epsilon}^{1}f(x)$ instead of $V_{t,\epsilon}f(x)$. Since by    (\ref{up-11})  we have
$$
|V_{t,\epsilon}^{1}f(x)|\leq Ct^{-1+\delta},
$$
the convergence  \eqref{conv_int}  for $V_{t,\eps}^1f(x)$ easily follows from  (\ref{conv_loc}).

Since  $|f|$ is bounded,   by \eqref{ld11} and  \eqref{Psi-est}  we obtain
\begin{equation}\label{int-sing}
\begin{split}
\int_\rn &\int_\rn \int_t^{t+\epsilon} |\Phi_{t-s+\epsilon} (x,z)\Psi_s(z,y) f(y)|dzdsdy  \\
&\leq c_1  \int_\rn\int_\rn  \int_t^{t+\eps} (t-s+\eps)^{-1+\eta} \big( g_{t-s+\eps}* G_{t-s+\eps} \big)(z-x) s^{-1+\delta} \big(g_{s,\zeta}*\Pi_s \big)(y-z)dsdzdy\\
&\leq  c_2 \int_t^{t+\eps} \int_\rn (t-s+\eps)^{-1+\eta}s^{-1+\delta} g_{t-\eps, \chi} (y-x-w) \Big[ \int_\rn
G_{t-s+\eps}(dw-u) \Pi_s(du)\Big] ds\\
&\leq c_3 \int_t^{t+\eps}  (t-s+\eps)^{-1+\eta}s^{-1+\delta}ds\\
&\leq c_4 t^{-1+\delta} \eps^{\eta},
\quad t\in [\tau,T], \quad x,y\in \rn.
\end{split}
\end{equation}
This immediately  gives  (\ref{conv_loc})  and (\ref{conv_int}) with   $V_{t,\epsilon}^{2}f(x)$ instead of $V_{t,\epsilon}f(x)$. \normal

\end{proof}

\subsection{Positive maximum principle, applied to the approximate fundamental solution. Proof of Theorem~\ref{t2}}\label{s52}

The proof of Theorem~\ref{t2} follows from Lemmas~\ref{posit}--\ref{Ptint} given below. The arguments used in the  proofs of these lemmas is literally the same as those used in  \cite{KK14a}. In order to make our paper self-contained, we give the proof of Lemma~\ref{posit}, a hint for the proofs of Lemmas~\ref{semigr} and \ref{Ptint}, and refer to \cite{KK14a} for the details.
\begin{lemma}\label{posit}
The operator  $S_t$ defined in \eqref{t-semi} is positivity preserving, i.e. $S_t f \geq 0 $ if $f\geq 0$.
\end{lemma}

\begin{proof}
Take $f\in C_\infty(\rn)$, $f\geq 0$, and suppose that
\begin{equation}\label{ass}
\inf_{t,x} S_t f(x)<0.
\end{equation}
Then there exists $T>0$ such that
$$
\inf_{t\leq T,x\in \rn} S_t f(x)<0.
$$
Then by \eqref{conv_pte} there  exist $\varsigma>0, \theta>0, \eps_1>0$ such that
$$
\inf_{t\leq T,x\in \rn} \Big(S_{t,\eps} f(x)+\theta t\Big)<-\varsigma, \quad \eps<\eps_1.
$$
Denote
$$
u_{\epsilon}(t,x)= S_{t,\eps} f(x)+\theta t,
$$
and note that by \eqref{conv_x}
$$
u_{\epsilon}(t,x)\to \theta t>0, \quad \|x\|\to\infty
$$
uniformly in $t\in [0, T]$. Hence the above infimum is in fact attained at some point; in what follows we fix one such a point for each $\eps$, and denote it by  $(t_\eps, x_\eps)$.

Since $f(x)\geq 0$, by statement 2 of Lemma \ref{aux-ep} there exist $\eps_0>0$, $\tau>0$ such that
$$
S_{t,\eps} f(x)+\theta t\geq -{\varsigma\over 2}, \quad t\leq \tau, \quad \eps<\eps_0, \quad x\in \rn.
$$
Because
$$
u_{\epsilon}(t_\eps,x_\eps)=\min_{t\in [0, T], x\in \rn}u_\eps(t,x)<-\varsigma<-{\varsigma\over 2},
$$
we have   $t_\eps>\tau$ as soon as $\eps<\eps_0$.

The operator $L$ satisfies the positive maximum principle; that is, if  whenever $f\in D(L)$, and  $f(x_0)\geq 0$ where $x_0=\arg\max f(x)$, then $Lf(x_0)\leq 0$, cf. \cite[Ch.~4.2]{EK86}. Therefore
$$
L_x u_\eps(t,x)|_{(t,x)=(t_\eps,x_\eps)}\geq 0.
$$
In addition, for $\eps<\eps_0$
$$
\prt_t u_\eps(t,x)|_{(t,x)=(t_\eps,x_\eps)}\leq 0.
$$
Note that inequality sign here may appear if $t_\eps=T$, and because we have excluded another ``boundary case'' $t_\eps=\tau$, the inverse inequality is impossible.

Then
\begin{equation}\label{PMP}
(\prt_t-L_x) u_\eps(t,x)|_{(t,x)=(t_\eps,x_\eps)}\leq 0.
\end{equation}
On the other hand, because $t_\eps\in [\tau, T]$, $\eps<\eps_0$ we have  by the first  statement of  Lemma \ref{l5}
$$
(\prt_t-L_x) u_\eps(t,x)|_{(t,x)=(t_\eps,x_\eps)}=\theta +V_{t_\eps,\eps} f(x_\eps)  \to \theta>0, \quad \eps\to 0.
$$
This gives contradiction and shows that \eqref{ass} fails.
\end{proof}

\begin{lemma}\label{semigr}
The family of operators possesses the semigroup property: $S_{t+s} = S_sS_t$.
\end{lemma}
\begin{proof}
 Take $f\in C_\infty(\rn)$. Applying the same argument as that,   used in the proof of  Lemma~\ref{posit}, to the functions
$$
u_{\pm}(t,x)=\pm S_{t+s}f(x)\mp S_tS_sf(x),
$$
one can  show that $u_\pm (t,x)\geq 0$, which implies  the  identity
$S_{t+s}f(x)-S_tS_sf(x)=0$.

\end{proof}

\begin{lemma}\label{Ptint}
We have
\begin{itemize}

  \item[a)]
 \begin{equation}\label{Dy1}
S_tf(x)-f(x)=\int_0^tS_sLf(x)\, ds, \quad f\in C_\infty^2(\rn);
\end{equation}

\item[b)]  $$ S_t 1 =1. $$
\end{itemize}

\end{lemma}
\begin{proof}
Applying  the  same argument as that used in the proof of Lemma~\ref{posit} to the functions
 $$
u_{\pm}(t,x)=\pm (S_tf(x)-f(x))\mp \int_0^tS_sLf(x)\, ds, \quad f\in C_\infty^2(\rn),
$$
and using the second statement of Lemma~\ref{l5},
we get identity \eqref{Dy1}.

Statement b) follows from a)  by taking
 $f_n\to 1$, $f_n\in C_\infty^2(\rn)$,  such that $Lf_n(x)\to 0$.
\end{proof}

\normal
\begin{proof}[Proof of  Theorem \ref{t2}.]

By Lemmas~\ref{posit}, \ref{semigr}, and  \ref{Ptint},  the family of operators $(S_t)_{t\geq 0}$ forms a strongly continuous contraction semigroup, which is positivity preserving.  Since the semigroup $(S_t)_{t\geq0}$ possesses the  continuous transition probability density $p_t(x,y)$, the respective Markov process $X$ is the strong Feller.  Finally, expression \eqref{Dy1}  and the second statement of Lemma~\ref{P-cont} implies that the restriction of the generator of $(S_t)_{t\geq 0}$ coincides with $L$ on functions from $C_\infty^2(\rn)$.

\end{proof}

\section{Time derivatives. Proof of Theorem~\ref{t5}}\label{pr-t5}

Proposition~\ref{p0-der} allows  us to transfer the differentiability properties of $p_t^0(x,y)$ to $p_t(x,y)$. For this we need to establish the  continuity  and upper estimates on $\partial_t \Phi^{\star k}$ and $\partial_t \Psi^{\star k}$.

\begin{lemma}\label{phi-der}
The function $\Psi_t(x,y)$ is differentiable with respect to $t$, $\partial_t \Psi_t(x,y)$ is  continuous in   $(t,x,y)\in(0,\infty)\times \rn\times \rn$, and for any $T>0$  there exists a family of   sub-probability measures  \,  $\{ \Theta_t, \, t\geq 0\}$,  such that
\begin{equation}\label{phi-der10}
\big|\partial_t \Psi_t(x,y)\big|\leq t^{-2+\delta}(g_{t,\zeta}* \Theta_t)(y-x), \quad   t\in (0,T] , \quad x,y\in \rn.
\end{equation}
\end{lemma}
\begin{proof}
The  proof  follows the same strategy  as that of Theorem~\ref{t-main1}.  Using Proposition~\ref{p0-der}, we can obtain the estimate for $\partial_t \Phi_t(x,y)$ in the same way as it was done for $\Phi_t(x,y)$ in Lemma~\ref{Phi-up}:
\begin{equation} \label{der-Phi1}
\big|  \partial_t \Phi_t(x,y)\big|\leq C t^{-2+\delta}   \big(\tilde{g}_t *\mathcal{G}_{t}\big)(y-x),\quad t\in (0,T], \, x,y\in \rn,
\end{equation}
where  $C>0$, $\delta\in (0,1)$ is the same as in \eqref{Phi-up20},  $\tilde{g}_t$ is of the form \eqref{g20}, and  the family of measures $\mathcal{G}_{t} (du)$ is given by
\begin{equation}\label{G11}
\mathcal{G}_t(du):=  c\Big( \mathcal{P}_{t,\kappa}(du)+ (\Lambda_t*
\mathcal{P}_{t,\kappa})(du)\Big),
\end{equation}
$$
\mathcal{P}_{t,\kappa}(du):= \big(1+\rho_t^\kappa ( \|u\|^\kappa \wedge 1)\big)  \mathcal{P}_t(du).
$$
Here  $c>0$ is the normalizing constant, such that    $\mathcal{G}_t(\rn) \leq 1$ for all $t\in (0,T]$.

Note that by definition
\begin{equation}\label{GG}
\mathcal{G}_t \geq G_t.
\end{equation}

To show that $\prt_t\Phi_t(x,y)$ is continuous in $(t,x,y)$, we follow line by line the proof of continuity of $\Phi_t(x,y)$ (cf. the proof of Theorem~\ref{t-main1}.b).  Observe that the function  $\prt_t p_t^0(x,y)$ is continuous in $(t,x,y)$, and $\prt_t p^0_t(\cdot,y)\in C_\infty^2(\rn)$. Then for any $t_0>0$
 $$
 \Big| \prt_t p_t^0(x+u,y)- \prt_t p_t^0(x,y)-\nabla_x \prt_t p^0_t(x,y) \cdot u \I_{\{\|u\|\leq 1\}} \Big|\leq C\big(\|u\|^2\wedge 1 \big), \quad t>t_0, \, x,y\in \rn
 $$
 Therefore, the function  $\prt_t L_x p_t^0(x,y)=  L_x \prt_t p_t^0(x,y)$ is continuous, which together with  continuity of $\prt_t^2 p_t^0(x,y)$ implies continuity of $\prt_t  \Phi_t(x,y)$ in $(t,x,y)$.

To show the continuity of $\prt_t \Phi_t^{\star k}(x,y)$ for $k\geq 2$ we use induction. \normal  Write
\begin{equation}\label{Phik}
\Phi^{\star (k+1)}_t(x,y)=\int_0^{t/2}\int_{\rn}\Phi_{t-s}^{\star k}(x,z)\Phi_s(z,y)\,dz ds+\int_0^{t/2}\int_{\rn}\Phi_{s}^{\star k}(x,z)\Phi_{t-s}(z,y)\,dz ds.
\end{equation}
 Observe that now the functions under the integrals do not have singularities in $t$.
 Differentiating the above expression in $t$  we get \normal
\begin{equation}\label{Phik10}
\begin{split}
\prt_t\Phi^{\star (k+1)}_t(x,y)&=\int_0^{t/2}\int_{\rn}(\prt_t\Phi^{\star k})_{t-s}(x,z)\Phi_s(z,y)\,dz ds\\
&\quad +
\int_0^{t/2}\int_{\rn}\Phi_{s}^{\star k}(x,z)(\prt_t\Phi)_{t-s}(z,y)\,dz ds\\
&\quad +\int_{\rn}\Phi_{t/2}^{\star k}(x,z)\Phi_{t/2}(z,y)\, dz.
\end{split}
\end{equation}
 Since by the induction assumption all functions under the integrals are continuous in $(t,x,y)$,  the above expression implies the continuity of $\prt_t \Phi_t^{(k+1)}(x,y)$  in $(t,x,y)\in (0,\infty)\times\rn\times \rn$.  \normal

Let us show by induction that
\begin{equation}\label{Phik20}
\big|   \partial_t \Phi_t^{\star k} (x,y)\big|\leq \tilde{C}_{k}    t^{-2+k\delta}
\big(g_{t}^{(k)} *\mathcal{G}_t^{(k)}\big)(y-x),\quad k\geq 2,
\end{equation}
where the sequence $g^{(k)}_t$ is given by \eqref{gt-new}, and
\begin{equation}\label{calgtk}
\begin{split}
\mathcal{G}^{(k)}_t(dw):&=\frac{1}{  1+B((k-1)\delta,\delta) }\Big(  \int_0^1 \int_\rn r^{-1+\delta}(1-r)^{-1+\delta(k-1)} \mathcal{G}_{t(1-r)}^{(k-1)}(dw-u) \mathcal{G}_{tr} (du)dr \\
&\quad +\big( \mathcal{G}_{t/2}^{(k-1)} * \mathcal{G}_{t/2}\big)(dw)\Big), \quad k\geq 2.
\end{split}
\end{equation}
Suppose that \eqref{Phik20} holds true for some $k\geq 2$.   Using \eqref{der-Phi1}, \eqref{Phik10}, \eqref{Phi10a} and Lemmas~\ref{phi12}, \ref{con-g-new}, we derive \begin{align*}
\big| \prt_t &\Phi_t^{\star (k+1)} (x,y)\big|  \leq \int_0^{t/2}\int_{\rn}\big|(\prt_t\Phi^{\star k})_{t-s}(x,z)\Phi_s(z,y)\big| \,dz ds\\
&\quad +
\int_0^{t/2}\int_{\rn}\big|\Phi_{s}^{\star k}(x,z)(\prt_t\Phi)_{t-s}(z,y)\big|\,dz ds +\int_{\rn}\big|\Phi_{t/2}^{\star k}(x,z)\Phi_{t/2}(z,y)\big|\, dz\\
&\leq c_1(k) \int_\rn g_t^{(k+1)} (y-x-w) \Big[ \int_0^{t/2} \int_\rn (t-s)^{-2+k\delta} s^{-1+\delta} \mathcal{G}_{t-s}^{(k)}(dw-u) G_s(du)ds\Big]\\
&+ c_2(k) \int_\rn g_t^{(k+1)} (y-x-w) \Big[ \int_0^{t/2} \int_\rn  (t-s)^{-1+k\delta} s^{-2+\delta}{G}_{t-s}^{(k)}(dw-u) \mathcal{G}_s(du)ds\Big]\\
&+ c_3(k) t^{-2+(k+1)\delta}\int_\rn g_t^{(k+1)}(y-x-w) \Big[ \int_\rn G_{t/2}^{(k)}(dw-u)G_{t/2}(du)\Big]\\
&\leq c_4(k) \int_\rn g_t^{(k+1)} (y-x-w) \Big[ \int_0^{t} \int_\rn (t-s)^{-1+k\delta} s^{-1+\delta} \mathcal{G}_{t-s}^{(k)}(dw-u) \mathcal{G}_s(du)ds\\
&\quad +t^{-2+(k+1)\delta} \Big( \mathcal{G}_{t/2}^{(k)}*\mathcal{G}_{t/2}\Big)(dw)\Big]\\
&\leq \tilde{C}_k t^{-2+(k+1)\delta}
\big(g_{t}^{(k+1)} *\mathcal{G}_t^{(k+1)}\big)(y-x),
\end{align*}
where we used \eqref{GG}.  This proves \eqref{Phik20}.

 Take as before $k_0:= \big[ \tfrac{n}{\alpha \delta}\big] +1$. Then applying induction and  \eqref{gtk0-new} (cf. Lemma~\ref{gtkl}) we get
\begin{equation}\label{Phik30}
\big| \partial_t \Phi_t^{\star (k_0+\ell)}(x,y) \big|\leq  \tilde{D}_\ell  t^{-2+ \delta (k_0+\ell)} \big(g_{t,\zeta} * \mathcal{G}_t^{(k_0+\ell)} \big) (y-x), \quad \ell\geq 1,
\end{equation}
where
$$
\tilde{D}_\ell := \frac{C(k_0) K^\ell \Gamma^{k_0+\ell}(\delta)}{\Gamma((k_0+\ell)\delta)}, \quad \ell\geq 1,
$$
and $C(k_0), K>0$ are some constants.

Finally, define
\begin{equation} \label{Pit}
\Theta_{t}(du):=  \frac{ \sum_{k=1}^{k_0} \tilde{C}_k T^{\delta(k-1)}\mathcal{G}_{t}^{(k)}(du)+ \sum_{\ell=1}^\infty \tilde{D}_\ell T^{\ell(k_0+\ell-1)} \mathcal{G}_t^{(k_0+\ell)}(du)}{\sum_{k=1}^{k_0} \tilde{C}_kT^{\delta(k-1)} + \sum_{\ell=1}^\infty \tilde{D}_\ell T^{\ell(k_0+\ell-1)} }, \quad t\in (0,T].\normal
 \end{equation}
Then $\Theta_t(\rn)\leq1$, $t\in (0,T]$, and thus  \eqref{phi-der10} follows from \eqref{Phik20} and \eqref{Phik30}.
\end{proof}

\begin{proof}[Proof of Theorem~\ref{t5}]
The proof of differentiability of $p_t(x,y)$ essentially follows from Proposition~\ref{p0-der} and  Lemma~\ref{phi-der}. Indeed, writing $p_t(x,y)$ in the form
\begin{equation}\label{pt-dec}
p_t(x,y)
=p_t^0(x,y)+\int_0^{t/2}\int_{\rn}p_{t-s}^0(x,z)\Psi_s(z,y)\,dzds+\int_{0}^{t/2}\int_{\rn}p_{s}^0(x,z)\Psi_{t-s}(z,y)\,dz\, ds,
\end{equation}
and applying the above lemmas we get
$$
 \big|\partial_t p_t(x,y)\big|\normal \leq C t^{-1} \big(g_{t,\chi} * \tilde{Q}_t\big) (y-x),\quad t\in (0,T], \, x,y\in\rn,
$$
where    $\chi\in (0,\zeta)$\normal, $\zeta$ is coming from \eqref{Phik30},
$$
\tilde{Q}_t(du):= c\big(\mathcal{P}_t(du)+ t^\delta\tilde{\mathcal{P}}_t(du)\big),
$$
is the a     sub-probability   measure (here   $c=c(T)>0$\normal\,  is the normalizing constant),
\begin{align*}
\tilde{\mathcal{P}}_t(du)&=
(\mathcal{P}_{t/2}* \Pi_{t/2})(du)+\int_0^{1/2} \int_\rn r^{-1+\delta} \mathcal{P}_{t(1-r)}(dw-u) \Pi_{tr} (du) dr\\
&\quad +\int_0^{1/2} \int_\rn  \Theta_{t(1-r)}(dw-u)P_{tr} (du)dr,
\end{align*}
 the measures $\Pi_t(du)$ and  $\Theta_t(du)$ are given, respectively, in  \eqref{Pi} \normal and \eqref{Pit}.
\end{proof}

We finish this section with a lemma, which plays an important role in the proof of Theorem~\ref{t3}.

\begin{lemma}\label{dptx}

\begin{enumerate}

  \item  For any  $f\in C_\infty(\rn)$,
$$
\|\prt_tS_{t,\eps} f-\prt_t S_t f\|_\infty \to0,  \quad \epsilon\to 0,
 $$
uniformly on compact subsets of $(0,\infty)$. Moreover, $\prt_tS_tf(x)=\int_\rn \prt_t p_{t}(x,y)f(y)dy$.

 \item
 $$
 \partial_t p_{t,\epsilon}(x,y)\to \partial_t p_t(x,y)\quad \text{as}\quad \epsilon\to 0,
 $$
   uniformly on compact subsets of $(0,\infty)\times \rn\times \rn$.

\end{enumerate}

\end{lemma}
The proof relies on the decomposition \eqref{pt-dec}, and the estimates on  $p^0, \Psi, \prt_t p^0, \prt_t\Psi$ obtained above; see the proof of \cite[Lem.6.4]{KK14a} for details.

\section{Proofs of Theorems~\ref{t3} and \ref{mart}}\label{pr-gen}

The proofs repeat literally the proofs  of the respective statements in \cite{KK14a}. In order to make this paper  self-contained, we sketch this proof below.

\begin{proof}[Proof of Theorem~\ref{t3}]
By Theorem~\ref{t2} we already know that
$(L, C_\infty^2 (\rn))$ is the restriction of $(A, D(A))$. Since $A$ is  closed, this yields that $(L, C_\infty^2 (\rn))$ is closable. Let us show that its closure coincides with $(A, D(A))$.

Take $f\in C_\infty(\rn)\cap D(A)$.
Fix $t>0$, and consider the functions  $S_tf$ and $S_{t,\eps} f$.  Since $f\in D(A)$, then $S_t  f\in D(A)$, and
 \begin{equation}\label{Pta}
 AS_t f=\partial_t S_tf.
 \end{equation}
 Recall that $S_{t,\eps}f (\cdot)\in C_\infty^2(\rn)\subset D(A)$, which implies
  $$
  AS_{t, \eps} f=L S_{t,\eps}f=\partial_t S_{t,\eps}f.
  $$
 Further, statement 2 in  Lemma \ref{aux-ep} together with statement 1 of Lemma~\ref{l5} and statement 1 of Lemma~\ref{dptx} implies
$$
LS_{t,\eps}f\to AS_tf\quad \text{in $C_\infty(\rn)$ as $\eps\to 0$},
$$
and thus  $S_tf $ belongs to the domain of the $C_\infty(\rn)$-closure of $(L, {C}_\infty^2 (\rn))$.   Consequently, this closure coincides with  $(A, D(A))$.

In addition, applying the same argument to the function $p_{t,\eps}(x,y)$ instead of $S_{t,\eps} f(x)$ and using  the second statement of Lemma~\ref{dptx}, we derive that $p_t(x,y)$ belongs to $D(A)$, and  is the fundamental solution to the Cauchy problem for $\partial_t-A$.
\end{proof}

\begin{proof}[Proof of Theorem~\ref{mart}]
Using the Markov property of $X$, we deduce from \eqref{Dy1}  and the  semigroup property for $p_t(x,y)$ the following: For given $f\in C_\infty^2(\rn)$, $t_2>t_1$, and $x\in\rn$, for any $m\geq 1$, $r_1, \dots r_m\in [0, t_1]$, and bounded measurable $G:(\rn)^m\to \rn$ the identity
$$
\Ee_x\left[f(X_{t_2})-f(X_{t_2})-\int_{t_1}^{t_2}h_f(X_s)\, ds\right]G(X_{r_1}, \dots, X_{r_m})=0
$$
holds true. Thus, for every $f\in C_\infty^2(\rn)$ the process
$$
M^f_t=f(X_t)-\int_{0}^{t}h_f(X_s)\, ds, \quad t\geq 0,
$$
is a $\Pp_x$-martingale for every $x\in \rn$; that is, $X$ is a solution to the martingale problem for $(L, C_\infty^2(\rn))$.

Note that the operator $(L,C_\infty^2(\rn))$ is dissipative, which follows from the positive maximum principle, see \cite[Lem.~4.2.1]{EK86}, or \cite[Lem.~4.5.2]{Ja01}. Since its closure equals to the generator $A$ of the $C_\infty(\rn)$-semigroup  $\{S_t, t\geq 0\}$,  for every $\lambda>0$ the range of the resolvent $(\lambda -L)^{-1}$  in $C_\infty(\rn)$ is dense. Hence the required  uniqueness of the solution to the martingale problem $(L,C_\infty^2(\rn))$ follows by \cite[Thm.~4.4.1]{EK86}.\normal
\end{proof}

\section{Upper and lower bounds: Proof of Theorem~\ref{t4}}\label{up-lo}

\emph{Upper bound.}    The upper bound is essentially contained in the proof of Theorem~\ref{t-main1}. Namely, we already obtained the upper estimate on $p^0\star \Psi$, see \eqref{up-11}.
Combining this estimate with the estimate  \eqref{der-est} (for $k=0$) for $p^0$, we derive the upper bound in \eqref{eqII} with
\begin{equation}\label{Q}
Q_t(du):= (1+\delta)^{-1} \Big(P_t(du) + t^\delta \tilde{\Pi}_t(du)\Big).
\end{equation}
Since $P_t$ is the probability measure and $\tilde{\Pi}_t$ is the sub-probability measure for $t\in [0,T]$,
$Q_t(\rn)\leq 1$ \,   for all $t\in [0,T]$.\normal

\emph{Lower bound.}   By \eqref{up-11} and the fact that $\tilde{\Pi}$ is the sub-probability measure  we get
\begin{equation}
\big|(p^0\star \Psi)_t(x,y)\big|\leq c_1 \rho_t^n t^{\delta}, \quad \quad \text{ $x,y\in \rn$, $t\in (0,T]$},\label{ZF1}
\end{equation}
which implies the upper bound $p_t(x,x) \leq c_2 \rho_t^n$ for all $x\in \rn$ and $t\in (0,T]$. Finally, using Proposition~\ref{aux2} and \eqref{ZF1}, we derive for $t$ small enough
\begin{align*}
p_t(x,y)&\geq p_t^0(x,y)-\big|(p^0\star \Psi)_t(x,y)\big|
\geq \rho_t^n  f_{low}(\|y-x\|\rho_t) - c_1 \rho_t^n t^{\delta}
\\&
\geq c_2  \rho_t^n  f_{low}(\|y-x\|\rho_t).
\end{align*}
  \qed

\section{Proof of Theorem~\ref{loc-time}}\label{appl}

\begin{proof}[Proof of Theorem~\ref{loc-time}]
\emph{1. Sufficiency.}  We use the upper bound, constructed in Theorem~\ref{t4}.

Fix  $\ell \in\Sss^n$, and define $\theta_t:= \inf\{ r:\,\, q^U(r\ell )\geq 1/t\}$. Note that  by \textbf{A1}  we have  $\theta_t\asymp \rho_t$  for all $t\in (0,1]$. For any $T\in (0,1]$ we  derive, making the change of variables $s=\theta_t$ and using \eqref{kul} in the integration by parts,
\begin{equation}\label{u2}
\begin{split}
\int_0^T\int_\rn &p_t(x,y)\varpi(dy)dt\leq c_0\int_0^T \int_\rn\int_\rn  \theta_t^n e^{-c \|x-y-w\|\theta_t} \varpi(dy) Q_t(dw) dt\\
&\leq c_1 \int_{\theta_T}^\infty \int_\rn\int_\rn \frac{s^{n-1} q^L(\ell s)}{(q^U(s))^2} e^{-c \|x-y-w\|s}\varpi(dy) \mathcal{Q}_s(dw)ds\\
&\leq  c_2 \int_{\theta_T}^\infty \int_\rn\int_\rn \frac{s^{n-1}}{q^*(s)} e^{-c \|x-y-w\|s}\varpi(dy) \mathcal{Q}_s(dw)ds\\
&= c_2\int_{\theta_T}^\infty \frac{s^{n-1}}{q^*(s)} \int_\rn \int_0^\infty \varpi\{ y:\, e^{-c \|x-y-w\|s}>r\}dr \mathcal{Q}_s(dw)ds\\
&= c_3  \int_{\theta_T}^\infty \frac{s^{n-1}}{q^*(s)} \int_\rn \int_0^\infty \varpi\{y:\, \|x-y-w\|\leq v/s\} e^{-c v} dv \mathcal{Q}_s(dw)ds\\
&\leq c_4 \int_{\theta_T}^\infty \frac{s^{n-1}}{q^*(s)} \int_0^\infty h(v/s) e^{-c v} dv ds\\
&= c_4 \int_0^\infty\Big[  \int_0^{1/\theta_T} \frac{h(sv)}{s^{n+1}q^*(1/s)} ds\Big]e^{-c v}dv,
\end{split}
\end{equation}
where $\mathcal{Q}_s(dw)$ is the image measure of $Q_t(dw)$ under the transformation $s=\theta_t$,
$$
h(r):= \sup_{x\in \rn} \varpi\{B(x,r)\},
$$
 and in the second line from below  we used  that    $\mathcal{Q}_s(\rn)\leq 1$ for all $s\in (\theta_T,\infty]$ .    Without loss of generality assume that $c=1$. Split
 \begin{align*}
 I(T):&=\int_0^\infty\Big[  \int_0^{1/\theta_T} \frac{h(sv)}{s^{n+1}q^*(1/s)} ds\Big]e^{-v}dv\\
 &= \int_0^1 \Big[\int_0^{1/\theta_T} \frac{h(sv)}{s^{n+1}q^*(1/s)} ds \Big] e^{-v}dv +
 \int_1^\infty  \Big[\int_0^{1/\theta_T} \frac{h(sv)}{s^{n+1}q^*(1/s)} ds \Big] e^{-v}dv\\
 & = : I_1(T)+I_2(T).
  \end{align*}
We show that under \eqref{d-suf} (respectively, \eqref{d-nes}) one has $I(T)\to 0$ as $T\to 0$ (respectively, $I(T)<\infty$).

 By monotonicity of $h(r)$ and \eqref{d-suf}  we have
 $$
 I_1(T)\leq \int_0^{1/\theta_T}   \frac{h(s)}{s^{n+1}q^*(1/s)} ds  \cdot \int_0^1 e^{-v}dv\to 0 \quad \text{as $T\to 0$}.
 $$
 Further, using monotonicity of $q^*$, we get
 \begin{align*}
 I_2(T)&\leq \int_1^\infty \Big[\int_0^{v/\theta_T}  \frac{h(u)}{u^{n+1}q^*(1/u)} du \Big] v^n e^{-v}dv\\
 &=\Big[ \int_1^\infty \int_0^{1/\theta_T} + \int_{1}^\infty \int_{1/\theta_T}^{v/\theta_T} \Big]
  \frac{h(u)}{u^{n+1}q^*(1/u)} du \Big] v^n e^{-v}dv=: I_{21}(T)+I_{22}(T).
  \end{align*}
  For $I_{21}(T)$ we have
  $$
  I_{21}(T)= \int_1^\infty v^n e^{-v}ds\cdot \int_0^{1/\theta_T} \frac{h(u)}{u^{n+1}q^*(1/u)} du \to 0, \quad T\to 0.
  $$
  Further,
  \begin{align*}
  I_{22}(T)&=\int_{1/\theta_T} ^\infty \Big[\int_{u\theta_T}^\infty v^n e^{-v}dv\Big]   \frac{h(u)}{u^{n+1}q^*(1/u)} du \\
 &\leq \int_0^\infty e^{-\epsilon u \theta_T} \Big[ \int_{u\theta_T}^\infty v^n e^{-(1-\epsilon)v} dv \Big]  \frac{h(u)}{u^{n+1}q^*(1/u)} du\\
 &\leq \int_0^\infty e^{-\epsilon u } \Big[ \int_{u\theta_T}^\infty v^n e^{-(1-\epsilon)v} dv \Big]  \frac{h(u)}{u^{n+1}q^*(1/u)} du.
 \end{align*}
 Since  by \eqref{d-suf}  the function $\phi(u):=  \frac{e^{-\epsilon u }   h(u)}{u^{n+1}q^*(1/u)}$ is integrable on $(0,\infty)$, we have by the theorem on continuity with respect to a parameter that $I_{22}(T)\to 0$ as $T\to 0$. Therefore,  under \eqref{d-suf} (resp., \eqref{d-nes}) we have  $I(T)\to 0$ as $T\to 0$ (resp., $I(T)<\infty$) and thus $\varpi\in S_K$ (resp., $\varpi\in S_D$).

\emph{2. Necessity. }
Using the lower bound for $p_t(x,y)$ and the inequality $(1-\|x\|s)_+\geq 2^{-1} \I_{\{2\|x\|s\leq 1\}}$, we  obtain
\begin{equation}
\begin{split}
\int_0^T\int_\rn  p_t(x,y) \,\varpi(dy)dt&\geq d_1\int_0^T \int_\rn \rho_t^n(1-d_2\|x-y\|\rho_t)_+\,\varpi(dy)dt\\
&\geq 2^{-1} d_1 \int_0^T \int_\rn \rho_t^n \I_{\{2 d_2 \|x-y\|\rho_t\leq 1\}}\varpi(dy) dt.\label{qu44}
\end{split}
\end{equation}
Without loss of generality assume  that $\delta=\delta(T):= 1/\theta_T\in (0,1)$, and that  $2d_2=1$.
Therefore, using \eqref{kul} and \textbf{A1} we derive
\begin{equation}\label{qu5}
\begin{split}
\int_0^T \rho_t^n \I_{ \|x\|\rho_t\leq 1}ds &\geq  c_1\int_{1/\delta}^{1/\|x\|} s^{n-1} \frac{q^L(\ell s)}{(q^U(\ell s))^2}ds
\\
&\geq \beta^{-1} c_1\int_{1/\delta}^{1/\|x\|} \frac{s^{n-1}}{q^*(s)}ds
\\&
=  \beta^{-1} c_1 \Big(U(\|x\|)-U(\delta)\Big) \I_{\{\|x\|\leq \delta)\}},
\end{split}
\end{equation}
where
$$
U(r):= \int_1^{1/r} \frac{s^{n-1}}{q^*(s)}ds, \quad r\in (0,1).
$$
Performing integration by parts, we derive
\begin{align*}
\int_{\|x-y\|\leq \delta} \big( U(\|x-y\|)-U(\delta)\big) \varpi(dy)&= \int_0^{U(0)-U(\delta)} \varpi\{ y: \, U(\|x-y\|)\geq r+U(\delta)\} dr\\
&=\int_{U(\delta)}^{U(0)}  \varpi\{ y:\, U(\|x-y\|)\geq r\}dr\\
&=\int_0^\delta  \frac{\varpi\{ y:\, \|x-y\|\leq s\}}{s^{n+1} q^*(1/s)} ds.
\end{align*}
Note that $\delta_T\to 0$ if and only if $T\to 0$. Thus, if $\varpi\in S_K$  (respectively, $\varpi\in S_D$) then \eqref{k-nes} (respectively, \eqref{d-nes}) holds true.
\end{proof}

\section*{Appendix A}

\begin{proof}[Proof of   Proposition~\ref{lem1}  ]
Clearly, for $n=1$ the statement holds true.  For $n\geq2$ we have
\begin{equation}\label{qr}
\begin{split}
\mu\{ u:\, \|u\|\geq r\}&\leq \sum_{i=1}^n \mu\{ u:\, |u_i|\geq r n^{-1/2} \} \\
&\leq n \max_{1\leq i\leq n} \mu\{ u:\,   |u_i|\geq r n^{-1/2}\}  \\
&\leq n \max_{1\leq i\leq n} q^U( \sqrt{n}r^{-1}  \ell_i) \\
&\leq  n^2 \max_{1\leq i\leq n} q^U (r^{-1}  \ell_i)\\
&\leq n^2 q^* (1/r),
\end{split}
\end{equation}
  where $\ell_i:=(0,\ldots, \underset{i}{1}, \ldots 0) \in \mathbb{S}^n$, and in the third line we used that $q^U(\xi c) \leq (c^2 \wedge 1) q^U(\xi)$ holds true for any $c>0$.  Thus, for $n\geq 2$ we have
\begin{equation}
\Lambda_t (\rn) = t \mu\{ u:\, \|u\|\geq 1/\rho_t\} \leq n^2 t q^*(\rho_t)=  n^2. \label{la-es1}
\end{equation}
\end{proof}

\begin{proof}[Proof of Proposition~\ref{cor11}]
Using \eqref{growth2}  and \eqref{qr}  we derive
\begin{equation}\label{est1}
\begin{split}
\int_{\rho_t \|u\|\geq 1}\Big( \| u\|^\lambda \wedge 1\Big) &\mu(du)= \int_{1/\rho_t \leq \|u\|\leq 1} \|u\|^\lambda \mu(du) + \int_{\|u\|\geq 1} \mu(du)\\
&\leq\int_{1/\rho_t\leq \|u\|\leq 1} \int_0^{\|u\|^\lambda} dr\mu(du)+c_1\\
& \leq \iint \I_{1/\rho_t\leq \|u\|\leq 1} \I_{0<r<\|u\|^\lambda} dr\mu(du)+c_1\\
&\leq \int_0^{1/\rho_t^\lambda} \mu\{u:\, \|u\|\geq 1/\rho_t\}dr+ \int_0^1
\mu\{u:\, \|u\|\geq r^{1/\lambda}\}dr+c_1\\
&\leq n^2 \rho_t^{-\lambda} q^*(\rho_t) + \lambda  \int_1^{\rho_t}\frac{\mu\{ u:\, \|u\|\geq 1/r\}}{r^{1+\lambda}}dr+  c_1\\
&\leq   n^2 t^{-1} \rho_t^{-\lambda}+ \lambda  n^2    \int_1^{\rho_t}   \frac{q^*(r )}{r^{1+\lambda}} dr + c_1, \quad t\in (0,T].
\end{split}
\end{equation}
 \normal Note that condition \textbf{A1} implies   for any $\ell\in \Sss^n$  the inequalities
\begin{equation}
 q^U(r\ell) \leq q^*(r) \leq \beta  q^U(r\ell).\label{qul}
\end{equation}
Let us estimate $I_\ell (r):= \int_1^r \frac{q^U(v\ell)}{v^{1+\lambda}}dv$, where the vector  $\ell\in \Sss^n$ is fixed, and $r>1$.   Note that for any $\ell\in \Sss^n$  the mapping $r\mapsto q^U(r\ell)$ is absolutely continuous, and for any  $0<r_1<r_2$, $\ell\in \Sss^n$, we have
\begin{equation}\label{kul}
q^U(r_2 \ell)-q^U(r_1 \ell)= \int_{r_1}^{r_2} \frac{2 q^L(v\ell)}{v}dv.
\end{equation}
Recall that
 $\lambda\in [0,\alpha)$.  Therefore, applying  \eqref{kul} we derive
\begin{align*}
I_\ell (r)\leq \beta \int_1^r \frac{q^L(v\ell )}{v^{1+\lambda}}dv = \frac{\beta}{2} \int_1^r \frac{1}{v^{\lambda}} dq^U(v\ell)
\leq  \frac{\beta}{2}  \left( \frac{q^U(r\ell)}{r^{\lambda}} +  \lambda I_\ell(r)\right),
\end{align*}
which  gives $I_\ell(r)\leq (\alpha-\lambda)^{-1} q^U(r\ell )/r^{\lambda}$. Applying again \eqref{qul} we get for any $\ell \in \Sss^n$
$$
  \int_1^r  \frac{q^*(v )}{v^{1+\lambda}} dv\leq \beta  I_\ell (r)  \leq \frac{2}{\alpha(\alpha-\lambda)} \frac{q^*(r)}{r^{\lambda}},\normal
$$
which together with the last line in \eqref{est1} finally gives
\begin{equation}\label{moment1}
\begin{split}
 \int_\rn   \big( \|u\|^\lambda \wedge 1 \big)   \Lambda_t(du)&\leq n^2 \rho_t^{-\lambda}+ c_2 t \rho_t^{-\lambda} q^*(\rho_t) + c_1t \leq c_3 \rho_t^{-\lambda}+ c_1 t \\
&\leq c_4\rho_t^{-\lambda},\quad t\in (0,T],\normal
\end{split}
\end{equation}
 where for the last inequality  we again used that $\lambda<\alpha$, and hence $t\rho_t^\lambda\leq c$, $t\in [0,T]$\normal.   This proves the statement of Proposition~\ref{cor11}.

\end{proof}
\begin{proof}[Proof of Proposition~\ref{ptkappa}]
From Proposition~\ref{cor11} we  have for any $T>0$
\begin{equation}\label{gam1}
\rho_t^\kappa \int_{\rn}  \Big( \|u\|^\kappa \wedge 1\Big)   \Lambda_t(du)\leq
C, \quad t\in [0,T].
\end{equation}
  Then by   \eqref{gam1} and   Proposition~\ref{lem1}    we have
\begin{align*}
\rho_t^{\kappa}& \int_\rn  (\|u\|^\kappa  \wedge 1)  \Lambda_t^{*2}(du)\leq \rho_t^\kappa \int_\rn \int_\rn (\|u\|^\kappa\wedge 1) \Lambda_t(du-w)\Lambda_t(dw)
\\&
 \leq 2^\kappa  \Big[ \rho_t^\kappa \int_\rn (\|u-w\|^\kappa\wedge 1) \Lambda_t(du-w)\Lambda_t(dw)+  \rho_t^\kappa \int_\rn (\|w\|^\kappa\wedge 1) \Lambda_t(du-w)\Lambda_t(dw)\Big]
\\&
\leq 2^\kappa  \Big[ \rho_t^\kappa \int_\rn (\|v\|^\kappa\wedge 1) \Lambda_t(dv)\Lambda_t(\rn)+  \rho_t^\kappa \int_\rn (\|w\|^\kappa\wedge 1) \Lambda_t(dw)\Lambda_t(\rn)\Big]\normal\\
&\leq  2^{\kappa+1}    n^2   C,\quad  t\in [0,T]\normal,
\end{align*}
where in the second line we applied  the inequality
\begin{equation}\label{ab}
(a+b)^\kappa \leq 2^{\kappa} (a^\kappa+ b^\kappa),\quad a,b\geq 0.
\end{equation}
   Let us check that
\begin{equation}
\rho_t^{\kappa} \int_\rn   (\|u\|^\kappa \wedge 1)
\Lambda_t^{*m}(du)\leq  C \big(2^{\kappa+1} n^2\normal\big)^{m-1}, \quad  m\geq 2,\quad t\in  [0,T]\normal. \label{ind}
\end{equation}
Indeed, by induction we have
\begin{align*}
\rho_t^{\kappa} \int_\rn   (\|u\|^\kappa \wedge 1)
\Lambda_t^{*m}(du)&= \rho_t^\kappa \int_\rn \int_\rn (\|u\|^\kappa\wedge 1) \Lambda_t^{*(m-1)}(du-w)\Lambda_t(dw)\\
&\leq 2^\kappa \Big[ \rho_t^\kappa \int_\rn (\|u-w\|^\kappa\wedge 1) \Lambda_t^{*(m-1)}(du-w) \Lambda_t(dw)\normal\\
&\quad \quad +  \rho_t^\kappa \int_\rn (\|w\|^\kappa\wedge 1) \Lambda_t^{*(m-1)}(du-w)  \Lambda_t(dw)\normal\Big]\\
&\leq 2^\kappa \Big[    C 2^{(m-2)(\kappa+1)}  n^{2(m-1)}  + Cn^{2(m-1)} \Big]\\
&\leq C \big(2^{\kappa+1}   n^2  \big)^{m-1}.
\end{align*}
Finally, by \eqref{ind} we have
\begin{equation}
\begin{split}
\rho_t^{\kappa} \int_{\rn}  \Big( \|u\|^\kappa \wedge 1\Big)   P_t(du)&   =   e^{-\Lambda_t(\rn)}   \sum_{m=1}^\infty  \frac{\rho_t^\kappa}{m!} \int_\rn (\|u\|^\kappa \wedge 1)
\Lambda_t^{*m}(du)
\leq \frac{C e^{2^{\kappa+1}    n^2  }}{2^{\kappa+1}   n^2  }, \quad  t\in[0,T]\normal.\label{gam2}
\end{split}
\end{equation}
\end{proof}

\begin{proof}[Proof of Proposition~\ref{lem1-1}]

Take an arbitrary $\theta\in (0,1)$. Using \eqref{ab} and  the inequality
$z^\kappa e^{-z} \leq c_1e^{-\theta z}$, $z\geq 0$,
where  $c_1>0$ is some constant, we derive
\begin{align*}
(\|x\|^{\kappa} \wedge 1)   f_t(x)
&\leq c_22^\kappa \rho_t^{-\kappa} \Big[ \int_\rn (\|\rho_t (x-w)\|\wedge \rho_t)^\kappa g_t(x-w) P_t(dw)
\\&
\quad +\int_\rn  g_t(x-w)( \|\rho_t w\|\wedge \rho_t)^\kappa P_t(dw)\Big]\\
& \leq c_22^\kappa \rho_t^{-\kappa} \Big[ c_3 \int_\rn g_{t,\theta} (x-w) P_t(dw)
\\&
\quad +  \int_\rn  g_t(x-w)( \|\rho_t w\|\wedge \rho_t)^\kappa P_t(dw)\Big]\\
&\leq   c_4  \rho_t^{-\kappa}  \int_\rn g_{t,\theta}(x-w)(1+ ( \|\rho_t w\|\wedge \rho_t)^{\kappa} )  P_t(dw)
\\&
=c_4 \rho_t^{-\kappa}  \big(g_{t,\theta} * P_{t,\kappa}\big)(x), \quad   t\in (0,T].\normal
\end{align*}

\end{proof}

\section*{Appendix B}

\begin{proof}[Proof of Proposition~\ref{H-bound-a}] 1. Using \eqref{der-est} we derive
\begin{align*}
\Big|\int_\rn p^0_t(x,y)f(y)dy|&\leq C \int_\rn\int_\rn g_t(y-x-w) |f(y)|P_t(dw)dz\\
&= C \int_\rn\int_\rn g_t(z-w)|f(x+z)| P_t(dw)dz.
\end{align*}
Then the right-hand  side follows by \eqref{int-est} and the dominated convergence theorem.

2. By definition of $p^0_t(x,y)$ (cf. \eqref{pto}) we have
\begin{align*}
\Big| \int_\rn p^0_t(x,y)f(y)dy-f(x)\Big|&\leq  \Big|\int_\rn
\mathfrak{p}_t^x (y-x)\big(f(y)-f(x)\big)dy\Big|\\
&\quad +\Big| \int_\rn \big(\mathfrak{p}_t^y(y-x)-\mathfrak{p}_t^x(y-x)\big) f(y)dy\Big|\\
&=I_1(t,x)+I_2(t,x).
\end{align*}
Fix $\eps>0$. Then since $f$ is continuous, we have $|f(x)-f(y)|<\eps$ as soon as $\|x-y\|\leq \delta$ for some $\delta=\delta(\eps,x)$. Then
\begin{align*}
I_1(t,x)&\leq  \left(\int_{\|x-y\|\leq \delta}+ \int_{\|x-y\|>\delta}\right)\mathfrak{p}_t^x (y-x)\big|f(y)-f(x)\big|dy\leq C_1\left(\eps + \int_{\|z\|>\delta} (g_t* P_t)(z)dz\right),
\end{align*}
where we used Proposition~\ref{aux1} (see also \eqref{gt10} and \eqref{ft}),  and that $f(x)$ is bounded.
Note that
\begin{align*}
I_{11}(t,x)&:=\int_{\|z\|>\delta } (g_t* P_t\big)(z)dz\leq C_2
\int_{\|u\|\geq \delta \rho_t} \int_\rn e^{-c\|u-\rho_t w\|} P_t(dw)du\\
&=C_2\int_{\|u\|\geq \delta \rho_t} \int_\rn e^{-c\|u-v\|} P_t^\sharp(dv)du,
\end{align*}
 where $P_t^\sharp (dv) $ is the measure, obtained from $P_t(dw)$ by the change of variables $\rho_t w=v$. Observe that $P_t^\sharp (dv) $ is a sub-probability measure.  Therefore, $\sup_x I_{11}(t,x)\to 0$ as $t\to 0$, which in turn implies  that $\lim_{t\to 0} \sup_x I_{1}(t,x)\leq \eps$.

Let us estimate $I_2(t,x)$. Since for $a,b>0$ we have $|e^{-a}-e^{-b}|\leq |a-b| e^{-(a\wedge b)}$, we get by the H\"older continuity of $m(x,u)$ (cf. the representation of $q(x,\xi)$)
 \begin{align*}
 \big|\mathfrak{p}_t^x(y-x)-\mathfrak{p}_t^y(y-x)\big|&=(2\pi)^{-n}  \Big| \int_\rn e^{-i\xi (y-x)} \Big( e^{-t q(x,\xi)}-e^{-t q(y,\xi)|}\Big) d\xi \Big| \\
 &\leq  c_1   (|y-x|^\gamma \wedge 1)\Big|\int_\rn t q^U(\xi) e^{- c t q^U (\xi)}\Big) d\xi \Big|\\
 &\leq c_2 |y-x|^\gamma \rho_t^n, \quad t\in (0,1], \quad x,y\in \rn.
 \end{align*}
 Take now $\varsigma> \frac{n}{n+\gamma}$. Then
 \begin{align*}
 I_2(t,x)&\leq \left(\int_{\|y-x\|\leq \rho_t^{-\varsigma} } + \int_{\|y-x\|> \rho_t^{-\varsigma} }\right)
  \big|\mathfrak{p}_t^x(y-x)-\mathfrak{p}_t^y(y-x)\big| |f(y)|dy\\
  &\leq C \Big(\rho_t^{n-(n+\gamma)\varsigma} + \int_{\|u\|\geq \delta \rho_t^{1-\varsigma}} \int_\rn e^{-c\|u-w\|} P_t^\sharp(dw)du\Big).
 \end{align*}
 By our choice of $\varsigma$, both terms tend to 0 as $t\to 0$, uniformly in $x$. Thus,
\begin{align*}
\lim_{t\to 0} \sup_x \Big| \int_\rn p^0_t(x,y)f(y)dy-f(x)\Big| <\eps.
\end{align*}
Since $\eps>0$ is arbitrary, this implies statement b).
\end{proof}

\textbf{Acknowledgement.}  The authors are very grateful to the anonymous referee for careful reading and very valuable remarks.   The authors thank K. Bogdan, N. Jacob,  R. Schilling and M. Z\"ahle  for inspiring discussions  and helpful remarks, and gratefully acknowledge the DFG Grant Schi~419/8-1.  The first-names author gratefully acknowledges  the Scholarship of the President of Ukraine for young scientists (2012-2014).

\end{document}